\documentclass[11pt]{article}
\usepackage[margin=1in]{geometry}
\usepackage{amsmath}
\usepackage{amssymb}
\usepackage{amsbsy}
\usepackage{amsthm}
\usepackage{amsfonts}
\usepackage{thmtools}
\usepackage{thm-restate}
\usepackage{url}
\usepackage{mathrsfs}
\usepackage{multirow}
\usepackage{booktabs}
\usepackage{chngpage}
\usepackage{subcaption}
\usepackage{float}
\usepackage{afterpage}
\usepackage{subfiles}
\usepackage[math]{cellspace}
\usepackage[font={small}]{caption}
\usepackage{tikz,pgfplots}
\pgfplotsset{compat=1.16}
\usetikzlibrary{shapes.geometric,arrows,positioning,calc}
\usetikzlibrary{patterns}
\usepackage[linesnumbered,ruled,noend]{algorithm2e}
\usepackage{graphicx}
\usepackage{chngcntr}
\usepackage{bbm}
\usepackage{bm}
\usepackage{enumitem}
\usepackage{setspace}
\usepackage{xcolor}
\usepackage{mathtools}
\usepackage[square,sort,comma,numbers]{natbib}
\usepackage[colorlinks=true,breaklinks=true,bookmarks=true,urlcolor=black,citecolor=black,linkcolor=black,bookmarksopen=false,draft=false,hypertexnames=false]{hyperref}
\usetikzlibrary{decorations.pathmorphing,backgrounds,fit,matrix,decorations.pathreplacing}
\tikzset{snake it/.style={decorate, decoration={snake, amplitude=.3mm, segment length=2mm}}}

\theoremstyle{remark}

\newtheorem*{claim*}{Claim}
\newtheorem*{remark*}{Remark}
\newtheorem{remark}{Remark}

\theoremstyle{plain}

\newtheorem{corollary}{Corollary}

\theoremstyle{definition}
\newtheorem{definition}{Definition}

\theoremstyle{plain}

\declaretheorem[name=Example, style=definition]{example}
\renewcommand\thmcontinues[1]{Cont.}

\makeatletter
\def\toclevel@definition{0}
\def\toclevel@remark{0}
\def\toclevel@lemma{0}
\def\toclevel@theorem{0}
\def\toclevel@corollary{0}
\def\toclevel@example{0}
\def\toclevel@proposition{0}
\makeatother

\newcommand{\pset}{\mathscr P}

\DeclareMathOperator{\conv}{conv}

\DeclarePairedDelimiter\iprod{\langle}{\rangle}

\newcommand*{\pfstart}{\begin{proof}}
\newcommand*{\pfend}{\end{proof}}
\newcommand{\nc}{\neg}
\newcommand{\ec}{\mathord{\sim}}
\newcommand{\id}{\mathrm{id}}

\tikzstyle{vertex} = [circle, minimum size=0.1cm, inner sep=0pt, draw=black, fill=black]
\tikzstyle{circ} = [circle, minimum width=0.5mm, inner sep=0pt,draw,fill]
\tikzstyle{hcirc} = [circle,minimum width=5mm, inner sep=0pt,draw]
\tikzstyle{bcirc} = [circle, minimum width=1.5mm, inner sep=0pt,draw,fill]
\tikzstyle{bhcirc} = [circle, minimum width=1.5mm, inner sep=0pt,draw]
\tikzstyle{ept} = [circle,minimum width=0mm, inner sep=0pt, white]
\tikzstyle{txt} = [text width=1.3cm,draw,rounded corners=3pt]
\tikzstyle{ncirc} = [circle,draw=black, inner sep=1pt, minimum width=4mm]

\definecolor{darkblue}{RGB}{0,0,139}

\newcommand{\rev}[1]{\textcolor{black}{#1}}

\newenvironment{revision}
{\color{black}}
{}

\allowdisplaybreaks

\begin{document}
\onehalfspacing

\title{Set System Approximation for Binary Integer Programs: Reformulations and Applications}
\author{Ningji Wei\\
\small Department of Industrial, Manufacturing \& Systems Engineering, Texas Tech University\\
\small \url{ningji.wei@ttu.edu}}
\date{}
\maketitle

\begin{abstract}
  Covering and elimination inequalities are central to combinatorial optimization, yet their role has largely been studied in problem-specific settings or via no-good cuts. This paper introduces a unified perspective that treats these inequalities as primitives for set system approximation in binary integer programs (BIPs). We show that arbitrary set systems admit tight inner and outer monotone approximations, exactly corresponding to covering and elimination inequalities. Building on this, we develop a toolkit that both \rev{identifies} classical structural correspondences (e.g., paths vs.\ cuts, spanning trees vs.\ cycles) and extends polyhedral tools from set covering to general BIPs, including facet conditions and \rev{strengthened valid inequalities}. We also propose new reformulation techniques for nonlinear and latent monotone systems, such as auxiliary-variable-free bilinear linearization, bimonotone cuts, and interval decompositions. 
  \rev{Computational experiments on bilinear supplier--region selection and distributionally robust network site selection show how these tools produce exact reformulations, strengthen root relaxations, and support efficient cut generation.}
  Overall, this unified view clarifies inner/outer approximation criteria, extends classical polyhedral analysis, and provides broadly applicable reformulation strategies for nonlinear BIPs.\ \\

\end{abstract}


\section{Introduction}
\label{sec:intro}
Valid inequalities are central to combinatorial optimization, not only for strengthening formulations but also for revealing structural properties of feasible sets. Among them, two recurring families are the \emph{covering} and \emph{elimination}
\footnote{
\rev{
  Although packing is often paired with covering, we use elimination to emphasize that these inequalities exclude forbidden structures, and to align with established terminology such as cycle elimination constraints.}
}
inequalities, defined for a structured subset $T$ (e.g., paths, trees, cycles) as
\[
\sum_{i \in T} x_i \geq 1 \quad \text{(Covering)}, \quad \sum_{i \in T} x_i \leq |T|-1 \quad \text{(Elimination)}.
\]
These inequalities have been widely recognized, derived, and often strengthened into stronger forms across a broad range of combinatorial optimization problems, including transportation and routing problems \citep{bellmore1968traveling,park2020unmanned,achuthan1996new,kara2004note,gao2023exact}, network design and survivability \citep{gunluk1999branch}, facility location \citep{laporte1986exact,church1974maximal}, sequencing and scheduling models \citep{clark2006multi}, and interdiction games \citep{israeli2002shortest,lozano2017backward,Church2004,smith2020survey,wei2025supervalid}.

Traditionally, such constraints are derived from problem-specific structures: for example, enforcing coverage over required sets leads to covering inequalities, while eliminating cycles in the traveling salesman problem (TSP) naturally yields elimination constraints. Subsequent frameworks such as Logic-Based Benders Decomposition (LBBD) \citep{hooker2003logic} and Combinatorial Benders Decomposition (CBD) \citep{codato2006combinatorial} have unified these ideas by interpreting such inequalities as \emph{no-good cuts} derived from infeasible solutions under monotone objective functions. In these frameworks, covering and elimination inequalities are generated iteratively to refine the feasible region by excluding solutions encountered during the solution process.

Despite these advances, several important questions regarding the broader role of these inequalities remain open in the context of general binary solution spaces:

\begin{itemize}
  \item When do covering and elimination inequalities yield valid outer or inner approximations of binary solution spaces, and under what conditions are these approximations exact?
  \item What general conditions ensure that such inequalities define facets of the convex hull of the feasible region?
  \item \rev{How can covering and elimination inequalities be extended to characterize nonlinear or latent monotone feasible regions?}
\end{itemize}

This paper addresses all three inquiries by introducing a new perspective that treats covering and elimination inequalities as primitive tools for \emph{binary solution space approximation}, rather than as post hoc feasibility cuts. This shift enables a unified method that (i) \rev{identifies} structural correspondences, such as $s$-$t$ paths versus $s$-$t$ cuts in flow problems \citep{ford1956maximal}, spanning trees versus edge cuts in min-cut problems \citep{karger2000minimum}, and edge cuts versus bipartite subgraphs in max-cut \citep{barahona1986cut,barahona1985facets}; (ii) extends classical tools from the set covering literature, including \rev{polyhedral analysis \citep{balas1972set,balas1989setii,sanchez1998set} and various lifting techniques \citep{wei2022integer,nobili1989facets,wei2025supervalid}}, to general binary integer programs (BIPs); and (iii) introduces new reformulation tools by identifying latent \emph{monotone systems}, i.e., solution spaces closed under subset or superset inclusion, thereby enabling more flexible use of these inequalities across diverse problem settings. In addition, we also introduce a compact set of algebraic properties of set system operators, extending the classical clutter--blocker relationship \citep{edmonds1970bottleneck,cornuejols2001combinatorial} to a richer family of operators. This serves as a convenient toolkit for structural reasoning and simplifies several of our results.

\subsection{Related Works}
The covering and elimination inequalities are closely related to two prominent strands of research. The first concerns classical constraint families such as set covering inequalities \citep{balas1972set} and subtour elimination constraints for TSP \citep{bellmore1968traveling}, which are widely used in mixed integer programming problems to strengthen formulations and improve computational efficiency. The second relates to no-good cuts and their generalizations within the frameworks of LBBD \citep{hooker2003logic} and CBD \citep{codato2006combinatorial}, where valid inequalities are iteratively derived from infeasible solutions to refine the master problem. Both offer insights that motivate the development of our set system approximation framework.

Set covering inequalities are widely employed in different contexts, such as interdiction games \citep{israeli2002shortest,wei2021integer,lozano2017backward}, vehicle routing \citep{park2020unmanned}, network design \citep{gunluk1999branch}, power grid optimization \citep{yao2007trilevel}, and \rev{facility} location problems \citep{church1974maximal}, often enhancing the branch-and-cut implementation framework. Moreover, these covering inequalities have been shown to possess strong facet properties in multiple problem settings \citep{wei2022integer,sassano1989facial,sanchez1998set}. Similarly, subtour elimination inequalities have become a critical reformulation component in transportation \citep{achuthan1996new,kara2004note,gao2023exact}, production scheduling \citep{clark2006multi}, location-routing problems \citep{laporte1986exact}, and have demonstrated significant impact in computational efficiency when properly strengthened and implemented \citep{desrochers1991improvements,smith1977computational}. 

In many applications, covering and elimination inequalities are derived via Benders feasibility cuts \citep{bnnobrs1962partitioning}, typically followed by strengthening steps to eliminate big-$M$ constants. This led to the development of \emph{combinatorial Benders decomposition} (CBD) \citep{codato2006combinatorial}, which generalizes Benders methods to mixed-integer problems with binary and continuous variables. A further generalization is offered by \emph{logic-based Benders decomposition} (LBBD) \citep{hooker2003logic}, which accommodates nonlinear and combinatorial subproblems by deducing new constraints from infeasible subproblem outcomes via logic inference. Due to its generality, LBBD has found wide applications in areas such as transportation, production, supply chain management, and telecommunications \citep{booth2016logic,li2022novel,mohamed2023two,botton2013benders}, with comprehensive discussions available in \cite{hooker2023logic,rahmaniani2017benders}. A notable feature of both LBBD and CBD is the use of \emph{no-good cuts}, which are typically defined for an encountered infeasible solution \( x \in \{0,1\}^n \) with index set \( I := \{i \in [n] \mid x_i = 1\} \) as
\[
  \sum_{i \in [n] \setminus I} x_i + \sum_{i \in I}(1 - x_i) \geq 1,
\]
ensuring that at least one variable in \( x \) must flip to eliminate this infeasible solution from the feasible set. When the objective function exhibits monotonicity, stronger variants, termed \emph{monotone cuts}, become valid. These take the form
\[
  \sum_{i \in [n] \setminus I} x_i \geq 1 \quad \text{or} \quad \sum_{i \in I} x_i \leq |I| - 1,
\]
depending on whether the objective function is increasing or decreasing in a minimization setting.

\begin{table}[tb]
\hbadness=10001
\centering
\small
\renewcommand{\arraystretch}{1.5}
\begin{tabular}{p{3cm}p{6cm}p{6cm}}
\toprule
\textbf{Aspect} & \textbf{LBBD / CBD / No-Good Cuts} & \textbf{Proposed Framework} \\
\midrule
Goal & Refine the feasible set of an optimization model using valid cuts & Approximate and characterize feasible structures of a given binary space \\
Constraint Validity & Constraints must be logically valid for the original problem & Constraints may be invalid (especially for inner approximation) but serve structural representation \\
\rev{Practical Usage} & \rev{Dynamically generates no-good or monotone cuts from individual infeasible solutions during the solution process} & \rev{Statically analyzes structures and uncovers latent monotone systems, such as bimonotone or interval systems, for approximation and reformulation}\\
\bottomrule
\end{tabular}
\caption{Comparison between classical methods (LBBD, CBD, and no-good cuts) and the proposed set system approximation framework.}\label{tab:comparison}
\end{table}

While the constraints in our framework resemble monotone cuts in form, the underlying focus is quite different. LBBD, CBD, and no-good cuts mainly serve to dynamically generate valid inequalities from individual infeasible solutions to guide solver convergence. By contrast, our framework focuses on the structural approximation of the binary solution space itself. Rather than operating on encountered solution vectors, it analyzes combinatorial set systems deduced from the entire family of \rev{infeasible solutions}. Depending on the intended approximation purposes---outer or inner---the resulting constraints need not be valid in the traditional logical sense; instead, they serve as structural primitives for representing or approximating the solution space via latent monotone subsystems, which may require additional reformulation or decomposition steps to uncover (e.g., bimonotone cuts or interval system decompositions in Section~\ref{sec:bipreform}). Table~\ref{tab:comparison} summarizes the key differences between the two approaches.

\rev{
Another related line of work is clutter--blocker theory, which reveals fundamental combinatorial interactions between minimal set systems and the sets that intersect them \citep{edmonds1970bottleneck,cornuejols2001combinatorial}. This perspective is closely connected to the set system operators studied in this work. Our framework builds on this connection by extending these structural interactions to arbitrary set systems and using them to derive tight monotone approximations. A detailed comparison is provided in Section~\ref{sec:cbtheory}.
}

\subsection{Contributions}
Our contributions, organized around the three motivating inquiries in the introduction, are summarized below.

\begin{itemize}
  \item \textbf{Set System Approximation.} 
    We establish the tightest inner and outer approximations of arbitrary set systems by monotone systems (Theorem~\ref{thm:muembd}), and prove their exact correspondence with covering and elimination inequalities (Theorem~\ref{thm:reform}). This resolves the \emph{first inquiry} by showing precisely when such inequalities yield outer/inner approximations of binary solution spaces, and when these approximations become exact. In doing so, we also develop a compact cut–cocut algebra \rev{(Section~\ref{sec:alg})} that extends the classical clutter–blocker framework \citep{edmonds1970bottleneck,cornuejols2001combinatorial} to a broader family of operators.




  \item \textbf{Unified Analytical Toolkit.} 
    Building on this result, we develop a unified toolkit for analyzing binary integer programs (BIPs). First, we introduce a systematic method for identifying structural correspondences \rev{(e.g., cycles and claws vs.\ paths, neighborhoods vs.\ dominating sets, trees vs.\ cycles)} directly from the framework (Section~\ref{sec:dualstruct}). \rev{Second, we extend classical polyhedral tools from set covering to general BIPs, including a general facet-defining condition (Section~\ref{sec:polyhedral}) and strengthened valid inequalities (Section~\ref{sec:ineq}).} This addresses the \emph{second inquiry}, demonstrating how existing combinatorial and polyhedral insights fit within the unified set-system perspective.  

  \item \textbf{New Reformulations.} We introduce several reformulation methods for nonlinear and latent monotone systems that, to the best of our knowledge, are new. These include: 
    (i) a linearization of bilinear terms $\iprod{x, Ry}$ that avoids auxiliary variables when $R$ is structured (Example~\ref{eg:bilinear}), in contrast to the standard McCormick linearization \citep{mccormick1976computability};
(ii) bimonotone cuts, obtained by extending monotone cuts via the flipping map (Example~\ref{eg:signmincut});
(iii) interval system decompositions, which provide reformulations for piecewise-affine objectives \rev{(Example~\ref{eg:piecewise})}.
\rev{The computational experiments evaluate these structures through bilinear supplier--region selection and distributionally robust site selection.}
Collectively, these tools address the \emph{third inquiry} and extend the scope of monotone reformulations beyond classical settings.
\end{itemize}

\rev{
  We use the following notation throughout.
  An empty feasible region has value \(+\infty\) for minimization and \(-\infty\) for maximization.
  $G=(V,E)$ denotes a connected simple undirected graph with vertex set $V$ and edge set $E$. For $v\in V$, $N[v]$ is the closed neighborhood of $v$; for $V'\subseteq V$, $G[V']$ is the subgraph induced by $V'$. All edge and vertex weights are nonnegative unless stated otherwise. Given $n \in \mathbb Z_{++}$, $[n]$ denotes $\{1,2,\dots, n\}$. For a matrix $R$ indexed by $[n] \times [m]$, $R_{IJ}$ denotes the submatrix with rows $I$ and columns $J$, $\sigma(I,J)$ its entry sum, and $\bar{I} = [n]\setminus I$, $\bar{J} = [m]\setminus J$ the complement index sets. $\iprod{x,y}$ denotes the inner product of two vectors of equal dimension. 
For any set \(\Xi\), the indicator function \(\mathbb I_\Xi(\xi)\) returns \(1\) if \(\xi\in\Xi\) and \(0\) otherwise.
}

The remainder of the paper is organized as follows. Section~\ref{sec:cgalgebra} develops the main set system approximation results for binary solution spaces. Section~\ref{sec:struct} demonstrates how these results recover structural correspondences in classical combinatorial problems and extend polyhedral tools from the set covering literature to general BIPs. Section~\ref{sec:bipreform} introduces two complementary reformulation strategies based on latent monotone system identification. 
\rev{Section~\ref{sec:case} presents two case studies on bilinear supplier--region selection and distributionally robust network site selection.}
Finally, Section~\ref{sec:conclusion} offers concluding remarks and a summary of our findings. 
Technical proofs are deferred to Appendix~\ref{sec:proofs} to streamline exposition.

\section{Set System Approximation}
\label{sec:cgalgebra}

\begin{figure}[tbp]
    \centering
    \includegraphics[width=1\textwidth]{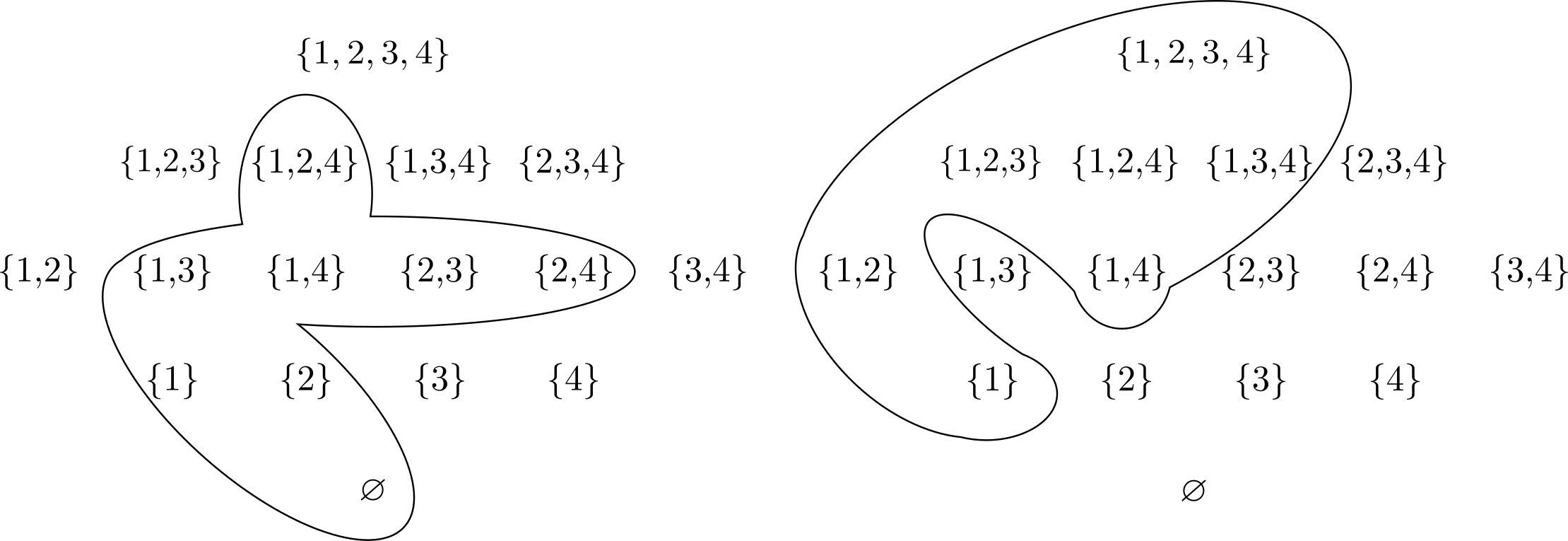}
    \caption{Two set systems, $\Omega_1$ (left) and $\Omega_2$ (right), each defined over the ground set $\Delta=\{1,2,3,4\}$. Each curve encloses all subsets belonging to $\Omega_1$ (left) or $\Omega_2$ (right).}
    \label{fig:poset}
\end{figure}

Every binary vector $x \in \{0,1\}^n$ can be uniquely associated with a subset of the ground set $\Delta := [n]$, defined by $T_x := \{i \in \Delta \mid x_i = 1 \}$. Conversely, for a subset $T \subseteq \Delta$, let $x_T$ denote the corresponding binary vector. Under this correspondence, any binary solution space $\mathcal X \subseteq \{0,1\}^n$ can be equivalently represented as a \emph{set system} $\Omega_{\mathcal X}:=\{T \subseteq \Delta \mid x_T \in \mathcal X\}$, which is a subset of the power set denoted as $\pset(\Delta)$. Figure~\ref{fig:poset} illustrates two set systems with $n=4$. 

This section develops the tightest inner and outer approximations of a given set system $\Omega$ using simpler structures, termed the \emph{monotone systems} defined as follows.

\begin{definition}[Monotone Systems]
A set system \( \Omega \subseteq \pset(\Delta) \) is \emph{upper-closed} if \( T \in \Omega \) and \( T' \supseteq T \) imply \( T' \in \Omega \); it is \emph{lower-closed} if \( T \in \Omega \) and \( T' \subseteq T \) imply \( T' \in \Omega \). Any upper- or lower-system is called a \emph{monotone system}.
\end{definition}

It is well-known that any union or intersection of upper-systems (or lower-systems) is still upper-closed (or lower-closed). To support the main approximation results, we introduce a set of basic operators that transform one set system into another.


\begin{definition} For a set system $\Omega$, we define eight operators organized into two groups.
  \label{defi:ops0}
  \begin{itemize}
    \item Closure \& Extremal Operators $\uparrow(\cdot), \downarrow(\cdot), m(\cdot), M(\cdot)$:
      \begin{itemize}
      \item  $\uparrow \Omega:=\{T  \mid T \supseteq T' \text{ for some } T' \in \Omega\}$ (\emph{Up-Closure Operator}),
      \item $\downarrow \Omega:=\{T  \mid T \subseteq T' \text{ for some } T' \in \Omega\}$ (\emph{Down-Closure Operator}),
      \item $m(\Omega):=\{T \in \Omega \mid \forall T' \in \Omega, T' \subseteq T \implies T' = T\}$ (\emph{Minimal Operator}),
      \item $M(\Omega):=\{T \in \Omega \mid \forall T' \in \Omega, T' \supseteq T \implies T' = T\}$ (\emph{Maximal Operator}).
      \end{itemize}
    \item Complement \& Cut Operators \rev{$\nc{(\cdot)}, \ec{(\cdot)}, \mathcal C(\cdot), \mathcal G(\cdot)$}:
      \begin{itemize}
       \item $\nc \Omega:= \pset(\Delta) \setminus \Omega$ (\emph{Complement Operator}),
       \item $\ec \Omega:= \{\Delta \setminus T \mid T \in \Omega\}$ (\emph{Element-Complement Operator}),
       \item $\mathcal C(\Omega) := \{S \mid \forall T \in \Omega, ~S \cap T \neq \emptyset\}$ (\emph{Cut Operator}),
       \item $\mathcal G(\Omega) := \{S \mid \forall T \in \Omega, ~S \cup T \neq \Delta\}$ (\emph{Cocut Operator}),
      \end{itemize}
  \end{itemize}
\end{definition}

\begin{revision}
\begin{remark}
 To simplify notation, we stack operators from right to left to denote they apply in this order. For instance $\nc \mathcal C m (\Omega)$ means first finding the minimal elements in $\Omega$, then apply the cut operator, followed by the complement. We also use $\id(\cdot)$ to denote the identity operator. Two operators are equal if they agree on every input.
\end{remark}
\end{revision}

The up- and down-closure operators return the smallest (with respect to inclusion) upper- or lower-systems containing $\Omega$, respectively. Conversely, the minimal and maximal operators remove redundant elements by retaining only the minimal or maximal sets in a monotone system. The complement operator returns all subsets of $\Delta$ not in $\Omega$ (i.e., structures associated with infeasible solutions), while the element-complement operator takes the complement of each $T \in \Omega$ within $\Delta$ (e.g., identifying all complement subgraphs when $\Delta$ is an edge set and $\Omega$ is a family of graphs defined by edge subsets). A \emph{cut} (also known as a \emph{hitting set} \citep{karp2009reducibility}) of a set system $\Omega$ is a subset of $\Delta$ that intersects every $T \in \Omega$. Then, the cut operator $\mathcal{C}(\Omega)$ collects the family of such sets and has seen applications in interdiction games \citep{wei2025supervalid}. 
This cut operator is closely related to the concepts of \emph{clutters} and \emph{blockers} \citep{edmonds1970bottleneck,cornuejols2001combinatorial}. \rev{A detailed comparison is provided in Section~\ref{sec:cbtheory}.}
In our setting, this cut operator will be used to derive approximations using upper systems. The cocut operator is the natural dual of the cut operator and will play a symmetric role in deriving approximations using lower systems. Its definition can be equivalently rewritten as:
$$
\begin{aligned}
  \mathcal{G}(\Omega) &= \{S \subseteq \Delta \mid \forall T \in \Omega, ~ (\Delta \setminus S) \cap (\Delta \setminus T) \neq \emptyset \}, \\
                      &= \{S \subseteq \Delta \mid \forall T \in \Omega, ~ \Delta \setminus S \not\subseteq T \}.
\end{aligned}
$$
The first form shows that $\mathcal{G}$ is the cut operator applied to the element-complement space $\ec{\Omega}$; the second reveals that each cocut in $\mathcal{G}(\Omega)$ ensures its complement is not fully contained in any structure $T \in \Omega$. 



\begin{revision}
\subsection{Cut-Cocut Algebra}
\label{sec:alg}
This subsection establishes basic algebraic properties of the eight set system operators. Although many of the results are intuitive, the resulting algebra provides a compact toolkit for streamlining later derivations and highlighting useful interactions among the operators. We begin by formalizing duality for set system operators.

\begin{definition}[Dual Operator]
Given a set system operator $g : \pset(\pset(\Delta)) \to \pset(\pset(\Delta))$, its \emph{dual operator} $g^\ast$ is defined by
\(
g^\ast = \ec g \ec
\).
\end{definition}

Intuitively, the dual operator \(g^\ast\) is obtained by conjugating \(g\) with the element-complement map: first transport \(\Omega\) to the element-complement space, apply \(g\), and then transport the result back. Equivalently, \(T\in g^\ast(\Omega)\) if and only if \(\Delta\setminus T\in g\ec(\Omega)\). Since the element-complement is an involution, i.e., \(\ec\ec=\id\), this operator duality is symmetric, i.e., $g^{\ast\ast} = \ec\ec g \ec\ec = g$.

A direct verification shows that \((\uparrow,\downarrow)\), \((m,M)\), and \((\mathcal C,\mathcal G)\) from Definition~\ref{defi:ops0} are dual pairs. The next lemma records the corresponding anti-commutativity relation with the element-complement operator, which will be used repeatedly. Its proof is omitted, as it follows directly from the definition of operator duality.

\begin{table}[t]
\hbadness=10001
\small
\centering
\rev{
\begin{tabular}{p{2.3cm}p{5cm}p{2.5cm}p{3cm}}
\toprule
\textbf{Operator} & \textbf{Basic Properties} & \textbf{Dual Operator} & \textbf{Interactions} \\
\midrule
Complement $\nc$ 
                  & Order-reversing; $\nc\nc=\id$; swaps upper/lower systems 
  & Self-dual & $\nc\ec=\ec\nc$ \\
\midrule
Element-complement $\ec$ 
  & Order-preserving; $\ec\ec=\id$; swaps upper/lower systems 
  & Self-dual & Anticommutes with dual operators \\
\midrule
Cut $\mathcal C$ 
  & Always upper; $\mathcal C\mathcal C=\uparrow$; depends only on minimal sets 
  & Cocut $\mathcal G$ 
  & $\mathcal C\ec=\ec{\mathcal G}$; $\ec\mathcal C={\mathcal G}\ec$ \\
\midrule
Cocut $\mathcal G$ 
  & Always lower; $\mathcal G\mathcal G=\downarrow$; depends only on maximal sets 
  & Cut $\mathcal C$ 
  & $\mathcal G\ec=\ec{\mathcal C}$; $\ec\mathcal G={\mathcal C}\ec$ \\
\bottomrule
\end{tabular}
}
\caption{Summary of set system operators and their algebraic properties.}
\label{tab:cheatsheet}
\end{table}

\begin{restatable}{lemma}{contra}
\label{lem:contra}
For any dual pair $g,g^\ast$, we have $\ec g = g^\ast \ec$ and $\ec g^\ast = g \ec$.
\end{restatable}

\begin{remark}
The notion of duality above is relative to the element-complement involution \(\ec\). More generally, any involution \(k\) on the space of set systems (e.g., $\nc$) induces a conjugate operator \(g^{k}:=kgk\), which satisfies the anti-commutativity relations
\(
kg = g^{k}k,
kg^{k}=gk
\) by the same argument.
Thus, operator duality can be viewed as conjugacy relative to a chosen involution. In this paper, we focus on $k:=\ec$ because it maps upper systems to lower systems and sends the cut operator to the cocut operator, thereby matching the covering--elimination symmetry used in the subsequent reformulations.
\end{remark}

Equipped with these preliminaries, we collect the basic properties for the complement, element–complement, cut, and cocut operators, along with their interactions in the following theorem. A summary is provided in Table~\ref{tab:cheatsheet}.


\begin{restatable}{theorem}{alg}
\label{thm:alg}
For any $\Omega \subseteq \pset(\Delta)$, the following hold.
\begin{enumerate}[label=\textbf{\Alph*}., leftmargin=2.2em]

\item \textbf{Complement $\nc$}
\begin{enumerate}[label=$\mathscr{A}\arabic*.$, ref=($\mathscr{A}\arabic*$)]
\item \label{alg:a1} $\nc{\emptyset}=\pset(\Delta)$ and $\nc(\pset(\Delta))=\emptyset$.
\item \label{alg:a3} $\Omega \subseteq \Omega'$ iff $\nc{\Omega} \supseteq \nc{\Omega'}$ (order-reversing).
\item \label{alg:a4} If $\Omega$ is upper (resp.\ lower), then $\nc{\Omega}$ is lower (resp.\ upper).
\item \label{alg:a5} $\nc{\nc}=\id$ (involution).
\end{enumerate}

\item \textbf{Element–complement $\ec$}
\begin{enumerate}[label=$\mathscr{B}\arabic*.$, ref=($\mathscr{B}\arabic*$)]
\item \label{alg:b1} $\ec{\emptyset}=\emptyset$ and $\ec(\pset(\Delta))=\pset(\Delta)$.
\item \label{alg:b2} $\Omega \subseteq \Omega'$ iff $\ec{\Omega} \subseteq \ec{\Omega'}$ (order-preserving).
\item \label{alg:b4} If $\Omega$ is upper (resp.\ lower), then $\ec{\Omega}$ is lower (resp.\ upper).
\item \label{alg:b5} $\ec{\ec}=\id$ (involution).
\end{enumerate}

\item \textbf{Cut $\mathcal C$}
\begin{enumerate}[label=$\mathscr{C}\arabic*.$, ref=($\mathscr{C}\arabic*$)]
\item \label{alg:c1} $\mathcal C(\emptyset)=\pset(\Delta)$ and $\mathcal C(\pset(\Delta))=\emptyset$.
\item \label{alg:c2} $\mathcal C(\Omega)=\emptyset$ iff $\emptyset\in\Omega$.
\item \label{alg:c3} $\mathcal C(\Omega)$ is an upper system.
\item \label{alg:c4} $\mathcal C=\mathcal Cm$.
\item \label{alg:c5} $\Omega \subseteq \Omega'$ implies $\mathcal C(\Omega)\supseteq \mathcal C(\Omega')$ (order-reversing).
\item \label{alg:c6} $\mathcal C\mathcal C=\uparrow$ (upper envelope).
\end{enumerate}

\item \textbf{Cocut $\mathcal G$}
\begin{enumerate}[label=$\mathscr{D}\arabic*.$, ref=($\mathscr{D}\arabic*$)]
\item \label{alg:d1} $\mathcal G(\emptyset)=\pset(\Delta)$ and $\mathcal G(\pset(\Delta))=\emptyset$.
\item \label{alg:d2} $\mathcal G(\Omega)=\emptyset$ iff $\Delta\in\Omega$.
\item \label{alg:d3} $\mathcal G(\Omega)$ is a lower system.
\item \label{alg:d4} $\mathcal G=\mathcal GM$.
\item \label{alg:d5} $\Omega \subseteq \Omega'$ implies $\mathcal G(\Omega)\supseteq \mathcal G(\Omega')$ (order-reversing).
\item \label{alg:d6} $\mathcal G \mathcal G=\downarrow $ (lower envelope).
\end{enumerate}

\item \textbf{Interactions.}
\begin{enumerate}[label=$\mathscr{E}\arabic*.$, ref=($\mathscr{E}\arabic*$)]
\item \label{alg:e1} $\ec{\nc}=\nc{\ec}$ (commutativity).
\item \label{alg:e3} $\ec\uparrow ={\downarrow}\ec$ and $\ec\downarrow ={\uparrow}\ec$ (anti-commutativity).
\item \label{alg:e2} $\ec{m} = M\ec$ and $\ec{M} = m\ec$ (anti-commutativity).
\item \label{alg:e4} $\ec\mathcal C = {\mathcal G}\ec$ and $\ec\mathcal G = {\mathcal C}\ec$ (anti-commutativity).
\item \label{alg:e5} $\nc{\uparrow (\Omega)} \subseteq \downarrow \nc{(\Omega)}$ and $\nc{\downarrow (\Omega)} \subseteq \uparrow \nc{(\Omega)}$ (partial anti-commutativity).
\item \label{alg:e6} $m\nc{(\Omega)}\subseteq \nc{M(\Omega)}$ and $M\nc{(\Omega)}\subseteq \nc{m(\Omega)}$ (partial anti-commute).
\item \label{alg:e7} $\mathcal G\nc{(\Omega)} \subseteq \nc{\mathcal C(\Omega)}$ and $\mathcal C\nc{(\Omega)} \subseteq \nc{\mathcal G(\Omega)}$, with equality iff $\Omega$ is upper and lower, respectively (partial anti-commutativity).
\end{enumerate}

\end{enumerate}
\end{restatable}

\end{revision}

\subsection{Tightest Monotone Approximations of Set Systems}
Building on the previously defined set operators, we present the main results for approximating a set system using upper-systems as follows.

\begin{restatable}[Upper Approximation]{theorem}{muembd}
  \label{thm:muembd}
 For every set system $\Omega$, we have
 $$\mathcal C\ec\nc(\Omega) \subseteq \Omega \subseteq \mathcal C\ec{\nc{(\uparrow\Omega)}},$$
where $\mathcal C\ec{\nc(\Omega)}$ and $\mathcal C\ec{\nc{(\uparrow\Omega)}}$ are the tightest inner and outer approximations of $\Omega$ using upper-systems, respectively. Moreover, equality holds throughout if and only if $\Omega$ is upper-closed.
\end{restatable}

This theorem shows that the two upper-systems, $\mathcal C\ec{\nc(\Omega)}$ and $\mathcal C\ec{\nc{(\uparrow\Omega)}}$, tightly sandwich $\Omega$, providing the best inner and outer approximations among all upper-systems. By symmetry, an analogous result holds for lower approximations as stated below. We omit its proof, as it follows directly from Theorem \ref{thm:muembd} by applying it to the element-complement space $\ec \Omega$.

\begin{corollary}[Lower Approximation]
  \label{coro:lembd}
 Given a set system $\Omega$, we have
 $$\mathcal G\ec{\nc(\Omega)} \subseteq \Omega \subseteq \mathcal G\ec{\nc{(\downarrow \Omega)}},$$
 where $\mathcal G\ec{\nc(\Omega)}$ and $\mathcal G\ec{\nc{(\downarrow\Omega)}}$ are the tightest inner and outer approximations of $\Omega$ using lower-systems, respectively. Moreover, equality holds throughout if and only if $\Omega$ is lower-closed.
\end{corollary}

Together, these results yield the tightest inner and outer approximations of any set system $\Omega$ using monotone systems. The following corollary extends these constructions to additional cases.

\begin{restatable}{corollary}{embdvar}
  \label{coro:embdvar}
  Given a set system $\Omega$, the tightest inner approximation by upper-systems with respect to $\nc\Omega$, $\ec\Omega$, and $\ec{\nc\Omega}$ are
$$\mathcal C\ec(\Omega) \subseteq \nc\Omega,~ \mathcal C\nc(\Omega) \subseteq \ec\Omega, ~\mathcal C(\Omega) \subseteq \ec{\nc\Omega},$$
and their equalities hold if and only if $\Omega$ is a lower system for the first two cases and is an upper \rev{system} for the third case. Symmetrically, $\mathcal G\ec(\Omega)$, $\mathcal G\nc(\Omega)$, and $\mathcal G(\Omega)$ are the respective tightest inner approximations by lower-systems with equality conditions reversed.
\end{restatable}

The following example illustrates these results by explicitly computing the inner and outer approximations for $\Omega_1$ and $\Omega_2$ in Figure~\ref{fig:poset}.

\begin{revision}
\begin{example}
Consider the set system \(\Omega:=\Omega_1\) in Figure~\ref{fig:poset}. We compute the upper-inner approximation \(\mathcal C\ec\nc(\Omega)\) step by step:
\[
\begin{aligned}
\nc(\Omega)
&= \big\{\{3\}, \{4\}, \{1,2\}, \{3,4\}, \{1,2,3\}, \{1,3,4\}, \{2,3,4\}, \{1,2,3,4\}\big\},\\
\ec\nc(\Omega)
&= \big\{\{1,2,4\}, \{1,2,3\}, \{3,4\}, \{1,2\}, \{4\}, \{2\}, \{1\}, \emptyset\big\},\\
\mathcal C\ec\nc(\Omega)
&= \emptyset,
\end{aligned}
\]
where the last equality follows from~\ref{alg:c2}, since \(\emptyset\in \ec\nc(\Omega)\). On the other hand, \(\uparrow\Omega=\mathcal P(\Delta)\), because \(\Omega\) contains the empty set. Hence,
\(
\ec\nc\uparrow(\Omega)=\emptyset\), and
\(
\mathcal C\ec\nc\uparrow(\Omega)=\mathcal P(\Delta)
\)
by~\ref{alg:c1}. Thus, for $\Omega_1$, both the tightest upper-inner and upper-outer approximations are trivial.
In contrast, the lower-inner approximation is
\[
\mathcal G\ec\nc(\Omega)=\{\emptyset,\{1\},\{2\}\}.
\]
Indeed, these are precisely the subsets whose union with every element of \(\ec\nc(\Omega)\) is still different from \(\Delta\). The lower-outer approximation is computed as
\[
\begin{aligned}
\downarrow\Omega
&= \big\{\emptyset, \{1\}, \{2\}, \{3\}, \{4\}, \{1,2\}, \{1,3\}, \{1,4\}, \{2,3\}, \{2,4\}, \{1,2,4\}\big\},\\
\nc\downarrow(\Omega)
&= \big\{\{3,4\}, \{1,2,3\}, \{1,3,4\}, \{2,3,4\}, \{1,2,3,4\}\big\},\\
\ec\nc\downarrow(\Omega)
&= \big\{\{1,2\}, \{4\}, \{2\}, \{1\}, \emptyset\big\},\\
\mathcal G\ec\nc\downarrow(\Omega)
&= \downarrow\big\{\{1,3\}, \{2,3\}, \{1,2,4\}\big\}.
\end{aligned}
\]
Therefore, \(\mathcal G\ec\nc(\Omega)\) and \(\mathcal G\ec\nc\downarrow(\Omega)\) are the tightest lower-inner and lower-outer approximations of \(\Omega\), respectively.
For the case \(\Omega:=\Omega_2\), a similar computation shows that the tightest upper-inner and upper-outer approximations are
\[
\mathcal C\ec\nc(\Omega)
= \uparrow\big\{\{1,2\},\{1,4\}\big\},
\qquad
\mathcal C\ec\nc\uparrow(\Omega)
= \uparrow\big\{\{1\}\big\},
\]
while the corresponding lower approximations are trivial.
\hfill\(\triangle\)
\end{example}
\end{revision}

\begin{revision}
\subsection{Connection to Clutter--Blocker Theory}
\label{sec:cbtheory}
The cut--cocut algebra and the associated set system approximation theorems are closely related to classical clutter--blocker theory, introduced by~\citet{edmonds1970bottleneck}.
At its core, a clutter \(\Omega\) is a set system consisting only of minimal elements under inclusion, and its blocker \(b(\Omega)\) is the minimal system that intersects every element of \(\Omega\). In the cut--cocut notation, this can be written as
\(
b = m\mathcal C,
\)
restricted to clutters. A fundamental identity established by~\citet{edmonds1970bottleneck} is
\(
bb(\Omega)=\Omega
\)
for every clutter \(\Omega\).
Beyond these connections, the cut--cocut algebra differs from the
clutter--blocker theory in two aspects.

First, they differ in scope. The clutter--blocker theory concerns only minimal set systems (clutters), whereas the cut--cocut algebra applies to arbitrary set systems. This extension enriches the study of set system operators in several ways: (i)~additional operators (Definition~\ref{defi:ops0}) are needed to complete the picture; (ii)~known identities admit finer decompositions: the classical $bb = \id$ on clutters refines, via properties~\ref{alg:c4} and~\ref{alg:c6}, to
$$
m\mathcal{C} m \mathcal{C} = m \mathcal{C}\mathcal{C} = m \uparrow = \id;
$$
(iii)~classical operators gain new interpretations, as
$\mathcal{C}\mathcal{C} = {\uparrow}$ exposes an envelope property
hidden in the blocker operator; 
and (iv)~the dual operator \(\mathcal{G}\), which arises naturally in our framework, is not usually treated as a standalone counterpart in classical clutter--blocker theory.

Second, the finer algebra makes the inverse problem tractable. The clutter--blocker literature addresses the forward question: given a clutter $\Omega$, what is its blocker? Theorem~\ref{thm:muembd} answers the inverse question: given an arbitrary set system $\Omega$, which systems $\Omega'$ satisfy $\mathcal{C}(\Omega') = \Omega$, or best approximate it? The complete algebra yields an inner solution $\ec\nc(\Omega)$ and an outer solution $\ec\nc{\uparrow}(\Omega)$, which coincide precisely when $\Omega$ is upper-closed. This inverse question falls outside the natural scope of the clutter--blocker
theory, whose operators act at a coarser level of granularity.

In this sense, the cut--cocut algebra is a natural extension of the clutter--blocker theory from minimal set systems to arbitrary ones.


\end{revision}

\subsection{Approximation-Based Reformulations for General BIPs}
The preceding approximation results offer two advantages: (i) they yield the tightest possible inner and outer approximations using upper and lower systems; and (ii) they are formulated using cut and cocut operators, which naturally correspond to classical covering and elimination inequalities, as shown in the following lemma.

\begin{restatable}{lemma}{ce}
  \label{lem:ce}
  Given any set system $\Omega$, the vector representations of $\mathcal C(\Omega)$ and $\ec{\mathcal C(\Omega)}$ correspond to the following \emph{covering} and \emph{elimination} inequalities, respectively.
  $$
  \begin{aligned}
    \mathcal X_{\mathcal C(\Omega)} &= \left\{x \in \{0,1\}^n~ \bigg\vert ~\sum_{i \in T}x_i \geq 1, \forall T \in \Omega\right\},\\
    \mathcal X_{\ec{\mathcal C(\Omega)}} &= \left\{x \in \{0,1\}^n~ \bigg\vert ~\sum_{i \in T}x_i \leq |T|-1, \forall T \in \Omega\right\}.
  \end{aligned}
  $$
Moreover, these two types of inequalities can be equivalently converted to each other by the substitution $y:= 1-x$.
\end{restatable}

\rev{We use \(+\infty\) and \(-\infty\) to represent infeasibility and unboundedness, respectively.
Then, the above lemma leads to the following approximation results for general binary integer programs (BIPs).}

\begin{restatable}{theorem}{reform}
  \label{thm:reform}
  Given a general BIP $z:=\min_{x \in \mathcal X}f(x)$, let $\Omega:=\{T_x \mid x \in \mathcal X \}$. The following reformulations serve as inner/outer approximations of the original problem,

\noindent\begin{minipage}{0.49\textwidth}
  \begin{subequations}
    \label{eq:pformub}
  {\setlength{\abovedisplayskip}{2pt}%
 \setlength{\belowdisplayskip}{2pt}%
  \begin{align} 
    \intertext{Upper-Inner Approximation $z_{ui}$:}
    \min\limits_{x \in \{0,1\}^n} &~ f(x)\\
    \text{s.t.} & ~\sum_{i \in T} x_i \geq 1, ~\forall T \in m\ec{\nc(\Omega)} \label{eq:pformub01}
\end{align}}
  \end{subequations}
  \end{minipage}
  \begin{minipage}{0.49\textwidth}
  \begin{subequations}
    \label{eq:pformuc}
  {\setlength{\abovedisplayskip}{2pt}%
 \setlength{\belowdisplayskip}{2pt}%
  \begin{align} 
    \intertext{Upper-Outer Approximation $z_{uo}$:}
    \min\limits_{x \in \{0,1\}^n} &~ f(x)\\
    \text{s.t.} & ~\sum_{i \in T} x_i \geq 1, ~\forall T \in m\ec{\nc{\uparrow(\Omega)}} \label{eq:pformuc01}
\end{align}}
  \end{subequations}
  \end{minipage}

\noindent\begin{minipage}{0.49\textwidth}
  \begin{subequations}
    \label{eq:pformlb}
  {\setlength{\abovedisplayskip}{2pt}%
 \setlength{\belowdisplayskip}{2pt}%
  \begin{align} 
    \intertext{Lower-Inner Approximation $z_{li}$:}
    \min\limits_{x \in \{0,1\}^n} &~ f(x)\\
    \text{s.t.} & ~\sum_{i \in T} x_i \leq |T|-1, ~\forall T \in m\nc(\Omega) \label{eq:pformlb01}
\end{align}}
  \end{subequations}
  \end{minipage}
  \begin{minipage}{0.49\textwidth}
  \begin{subequations}
    \label{eq:pformlc}
  {\setlength{\abovedisplayskip}{2pt}%
 \setlength{\belowdisplayskip}{2pt}%
  \begin{align} 
    \intertext{Lower-Outer Approximation $z_{lo}$:}
    \min\limits_{x \in \{0,1\}^n} &~ f(x)\\
    \text{s.t.} & ~\sum_{i \in T} x_i \leq |T|-1, ~\forall T \in m\nc{\downarrow(\Omega)} \label{eq:pformlc01}
\end{align}}
  \end{subequations}
  \end{minipage}
\ \\
with objective values satisfying
$z_{uo} \leq z \leq z_{ui}$ and $z_{lo} \leq z \leq z_{li}.$
  Moreover, \eqref{eq:pformub} and \eqref{eq:pformuc} are both equivalent to the original problem if and only if $\Omega$ is upper-closed; \eqref{eq:pformlb} and \eqref{eq:pformlc} are both equivalent to the original problem if and only if $\Omega$ is lower-closed.
\end{restatable}

In other words, 
\rev{every BIP admits both inner and outer approximations by pure covering or elimination models, yielding corresponding upper and lower bounds on the optimal value.}
These approximations are exact precisely when the feasible region itself is monotone, explaining why classical covering and elimination formulations succeed in problems like set covering and max-cut, but only approximate more general cases. 

\begin{revision}
Given a set system \(\Omega\) on ground set \(\Delta\), a membership oracle for \(\Omega\) is an algorithm that takes a subset \(T\subseteq\Delta\) as input and returns true if and only if \(T\in\Omega\). Depending on the problem, such an oracle may be implemented through a linear program, a discrete optimization problem, or a tailored algorithm. Let \(\tau(\Omega)\) denote the worst-case running time of one oracle call. The following corollary characterizes the integer point separation complexity: given any \(T\subseteq\Delta\), it bounds the time needed to either certify \(T\in\Omega\) or identify a violated inequality indexed by an element of \(m\ec{\nc\Omega}\) or \(m\nc\Omega\).

\begin{restatable}{corollary}{sep}
  \label{coro:sep}
Suppose \(\Omega \subseteq 2^\Delta\) is monotone with a membership oracle of complexity \(\tau(\Omega)\). Given any \(T \subseteq \Delta\), the integer point separation complexity for \eqref{eq:pformub}--\eqref{eq:pformlc} is of order \(O(|\Delta|\tau(\Omega))\).
\end{restatable}
\end{revision}

Building on these results, the next example introduces (to the best of the author's knowledge) a new linearization technique for bilinear terms $\iprod{x, Ry}$ with structured $R$. In contrast to the standard McCormick approach, which requires auxiliary variables, our reformulation avoids variable lifting altogether and instead leverages covering and elimination inequalities obtained from the monotone system perspective. 
\rev{
The goal is not to compete with McCormick linearization in compactness or generality, but to provide an exact alternative reformulation based on the underlying set system structure.
}
In Section~\ref{sec:polyhedral}, we further establish facet-defining conditions for these inequalities, and in \rev{Section~\ref{sec:case1} demonstrate their computational utility.}

\begin{example}[Linearization of Bilinear Terms]
  \label{eg:bilinear}
  For each supplier $i \in [n]$ and each region $j \in [m]$, $R_{ij} \geq 0$ represents the partnership value (e.g., number of end customers reachable, service compatibility, or market value unlocked) if supplier $i$ supports region $j$. Then, the following problem seeks to guarantee a partnership value level $\alpha$ while minimizing the total selection cost,
  $$
  \begin{aligned}
    \min_{x \in \{0,1\}^n, y \in \{0,1\}^m} &~ \iprod{c_1, x} + \iprod{c_2, y}\\
    \text{s.t.} &~ \iprod{x, Ry} \geq \alpha,\\
                &~ h(x, y) \geq 0,
  \end{aligned}
  $$
  where the last constraint encodes additional business requirements. A common reformulation method to treat the bilinear term $\iprod{x, Ry}$ is to introduce auxiliary variables $z_{ij} = x_iy_j$ to initialize the standard linearization step. 
  \rev{Since \(R\geq 0\), the feasible structures induced by \(\iprod{x,Ry}\geq \alpha\) form an upper system. Therefore, Theorem~\ref{thm:reform} yields the following exact reformulation:}
  \begin{revision}
  $$
  \begin{aligned}
    \min_{x \in \{0,1\}^n, y \in \{0,1\}^m} &~ \iprod{c_1, x} + \iprod{c_2, y}\\
    \text{s.t.} &~ \sum_{i \in I} x_i + \sum_{j  \in  J} y_j \geq 1, \quad \forall (I,J) \in m\ec\nc(\Omega)\\
                &~ h(x, y) \geq 0,
  \end{aligned}
  $$
where the structure set
\(
\Omega:=\left\{(I,J)\ \middle|\ \sum_{(i,j)\in I\times J} R_{ij}\geq \alpha\right\}
\)
indexes the feasible submatrices. 
Accordingly, \(\ec\nc\Omega\) consists of the row and column index sets $(I,J)$ whose removal leaves the submatrix $R_{\bar I, \bar J}$ with its total entry sum, denoted by $\sigma(\bar I, \bar J)$,  below \(\alpha\).
  \end{revision}
Since membership in $\Omega$ can be decided in polynomial time, this formulation can be solved via a cut-generation algorithm, using the efficient integer point separation procedure guaranteed by Corollary~\ref{coro:sep}.

Suppose the bilinear term appears in the objective function:
$$
\begin{aligned}
  \min_{x \in \{0,1\}^n,\, y \in \{0,1\}^m} &~ \iprod{x, Ry} \\
  &~ h(x,y) \geq 0.
\end{aligned}
$$
Let $z^\star$ be the optimal value and define $\Omega_{z^\star} := \{(I,J) | \sum_{(i,j) \in I \times J} R_{ij} \leq z^\star\}$. Since $R \geq 0$, this system is lower-closed. The problem can then be reformulated as
$$
\begin{aligned}
  \min_{x \in \{0,1\}^n,\, y \in \{0,1\}^m} &~ \iprod{R1, x} + \iprod{R^\intercal 1, y} \\
  \text{s.t.}\quad 
                                            & \sum_{i \in I} x_i + \sum_{j \in J}y_j \leq |I| + |J| - 1, 
      \quad \forall (I,J) \in \nc{\Omega_{z^\star}}, \\
    & h(x,y) \geq 0.
\end{aligned}
$$
This reformulation turns the problem into a feasibility problem: every feasible solution must already be optimal because of the elimination constraints, so the objective function is no longer essential for identifying an optimal solution.  
Although the exact value \(z^\star\) is unknown, it can be approached through a cut-generation procedure, \rev{termed the \emph{supervalid inequality} method in interdiction games~\citep{israeli2002shortest,wei2025supervalid}. These inequalities aggressively remove feasible solutions while retaining at least one optimal solution, either in the remaining feasible region or as an incumbent, and have shown computational advantages in several interdiction settings~\citep{salmeron2004analysis,kheirkhah2016improved,baycik2019robust}.}
Specifically, we relax the elimination constraints to form a master problem. Solving this master problem at iteration $k$ yields a feasible solution $(x_k,y_k)$ to the original problem with value $\bar z_k = \iprod{x_k,Ry_k}$. Let $\bar z:=\min_{k} \bar z_k$ be the value of the incumbent, we then assume that the true optimum must be strictly better, set $z^\star := \bar z - \epsilon$ for some sufficiently small $\epsilon > 0$, and generate one or more violated elimination constraints to exclude encountered solutions. This process iterates until the master problem becomes infeasible, at which point the best incumbent is guaranteed to be $\epsilon$-optimal. In this scheme, the proposed objective function serves primarily as a heuristic guide for producing high-quality candidate solutions in the relaxed master problem.
\hfill $\triangle$
\end{example}

This example demonstrates how Theorem~\ref{thm:reform} can be directly applied to reformulate bilinear terms \rev{when the associated feasible structures are monotone}, both in constraints and objectives, without introducing auxiliary variables. The resulting reformulation is complementary to existing approaches; for instance, its inequalities can be incorporated into other formulations to further strengthen their relaxations \rev{as demonstrated in Section~\ref{sec:case1}}.


More broadly, Theorem~\ref{thm:reform} yields several key implications that guide the rest of the paper. (i) It provides a unified analytical perspective to identify relevant structures, i.e., the pair $(\Omega, m\ec{\nc \Omega})$ for upper-systems and the pair $(\Omega, m\nc \Omega)$ for lower-systems. (ii) It provides a pathway to extend classical results from set covering problems to general BIPs. (iii) It promotes monotone systems as a tool for describing solution spaces, since each such system can be exactly characterized using covering or elimination inequalities, leading to reformulation methods to uncover latent monotone systems. We will explore (i) and (ii) in the next section and develop (iii) in Section~\ref{sec:bipreform}.

\section{Combinatorial and Polyhedral Structures in BIPs}
\label{sec:struct}
The approximation results developed in Section~\ref{sec:cgalgebra} provide a unified analytical framework for a broad class of binary integer programs (BIPs). This section illustrates its utility from two perspectives: 
\rev{(i) it provides a systematic procedure for identifying structural correspondences, thereby recovering many relationships that were previously derived in problem-specific contexts;}
and (ii) it extends classical polyhedral tools, particularly those from the set covering literature, to more general BIPs.

\subsection{\rev{Identifying Structural Correspondences}}
\label{sec:dualstruct}
\begin{table}[t]
\small
\centering
\renewcommand{\arraystretch}{1.2}
\resizebox{\textwidth}{!}{%
\begin{tabular}{ccc}
\toprule
\textbf{Problem} & \textbf{Feasible Set $\Omega$} & \textbf{Structural Counterpart}\\
                 && $m\ec{\nc{\uparrow(\Omega)}} \text{ for minimization; } m{\nc{\downarrow(\Omega)}} \text{ for maximization}$\\
\midrule
longest path & $s$-$t$ paths & cycles, claws, subgraphs preventing extension to $s$-$t$ paths\\
min dominating set & dominating sets & minimal closed neighborhoods\\
min spanning tree & spanning trees & global edge cuts\\
max spanning tree & spanning trees & simple cycles\\
\rev{sigmoid constraint} & \rev{feasible activation supports} & \rev{maximal response-deficient supports}\\

\bottomrule
\end{tabular}
}
\caption{Structural correspondences computed by the proposed framework. Each problem’s feasible set system $\Omega$ is paired with its structural counterpart: $m\ec{\nc{\uparrow(\Omega)}}$ for minimization or $m\nc{\downarrow(\Omega)}$ for maximization.}\label{tab:dualstruct}
\end{table}

\rev{When the underlying solution space is monotone}, Theorem~\ref{thm:reform} implies a strong connection between the structure pairs \( (\Omega, \ec{\nc\Omega}) \) and \( (\Omega, \nc\Omega) \), corresponding to the respective minimization and maximization problems. 
In particular, solving an optimization problem over one structural set requires imposing covering or elimination constraints indexed by the other.
\rev{In what follows, we illustrate these structural correspondences through examples drawn from graph-theoretic models and a nonlinear binary feasibility system. The key point is that these correspondences are obtained by the same operation, namely computing either \(\ec{\nc\Omega}\) or \(\nc\Omega\), depending on whether \(\Omega\) is upper- or lower-closed.}
 A summary of these structural relationships is provided in Table~\ref{tab:dualstruct}.

\begin{example}[Longest Path]
  \label{ex:lp}
  In the longest path problem with nonnegative edge lengths, the set system $\downarrow\Omega$ consists of all subgraphs that can be extended to an $s$-$t$ path. The corresponding structure set $\nc{\downarrow\Omega}$ includes subgraphs that cannot be extended to such a path. While this set lacks a simple and complete characterization, several sufficient conditions are readily verifiable. 
 \rev{In particular, several types of minimal subgraphs belong to $\nc{\downarrow\Omega}$, including simple cycles, claws, and subgraphs with more than one edge incident to either $s$ or $t$.
The first type of subgraph has been used to derive cycle- and clique-elimination cuts that improve computational performance~\citep{marzo2022new}. The claw structure has also appeared in the complexity analysis of longest path problems~\citep{li2000longest,broersma2013exact}. Although these structures do not fully describe \(\nc{\downarrow\Omega}\), they can be generated heuristically to strengthen relaxations for the associated binary integer program.
}
\hfill $\triangle$
\end{example}


\begin{example}[Minimum Dominating Set]
  A dominating set is a vertex subset $T \subseteq V$ such that every vertex in $V$ is either in $T$ or adjacent to some vertex in $T$. By definition, the family of dominating sets $\Omega$ is upper-closed. The associated structural system $m\ec{\nc\Omega}$ can be derived as:
  \begin{revision}
  \begin{itemize}
    \item $\nc\Omega$ consists of $T \subseteq V$ for which there exists a vertex $v \in V \setminus T$ that has no neighbors in $T$.
    \item $\ec{\nc\Omega}$ contains subsets $T$ that satisfies $T \supseteq N[v]$ for some $v \in T$.
    \item The minimal elements are the minimal closed neighborhoods: $m\ec{\nc\Omega} = m\{N[v]\mid v \in V\}$.
  \end{itemize}
  \end{revision}
  Thus, covering all closed neighborhoods is necessary and sufficient for vertex domination, a structural correspondence that has been previously explored in \cite{haynes2013fundamentals}.
\hfill $\triangle$
\end{example}

\begin{example}[Spanning Trees]
  In the spanning tree minimization problem, the upper system $\uparrow\Omega$ consists of all connected subgraphs. The associated structural set $m\ec{\nc{\uparrow\Omega}}$ then corresponds to all edge cuts---minimal sets of edges whose removal disconnects the graph. In contrast, when maximizing over spanning trees, we consider the lower system $\downarrow\Omega$, which consists of all forests (acyclic subgraphs). The corresponding structure set $m\nc\downarrow\Omega$ contains all simple cycles, as these are the minimal subgraphs that violate acyclicity. The first correspondence underlies the design of efficient min-cut algorithms based on spanning tree packings \citep{karger2000minimum}, while the second reflects the classical matroidal duality between bases (maximal independent sets) and circuits (minimal dependent sets) \citep{oxley2006matroid}.
\hfill $\triangle$
\end{example}

\begin{revision}
The preceding examples show how the framework identifies structural correspondences in familiar combinatorial models. We next turn to nonlinear binary constraints, where covering and elimination characterizations are less apparent, yet can still describe the feasible space exactly.

\begin{example}[Sum of Sigmoids Constraint]
Sigmoid functions, such as the logistic function \(\sigma(t)=e^t/(1+e^t)\), arise naturally in statistical models, utility functions, marketing, and machine learning, and have been incorporated into various optimization settings~\citep{udell2013maximizing,srivastava2014knapsack}. We consider a BIP setting with constraint learning~\citep{maragno2025mixed}, where a one-hidden-layer sigmoid model
\[
  f(x;\hat w,\hat\alpha)
  =
  \sum_{k\in[K]} \hat w_k
  \sigma\left(\hat\alpha_{k0}+\sum_{i\in[n]}\hat\alpha_{ki}x_i\right)
\]
predicts a performance level from a binary resource-activation vector \(x\in\{0,1\}^n\). We assume that the trained model is monotone in $x$, that is, \(\hat w_k\ge 0\) and \(\hat\alpha_{ki}\ge 0\) for all \(k\in[K]\) and \(i\in[n]\), while the intercepts \(\hat\alpha_{k0}\) may be unrestricted. The nonlinear constraint
\begin{equation}
  \label{eq:sigmoid}
f(x;\hat w,\hat\alpha)\ge \beta
\end{equation}
then requires the predicted performance to exceed a prescribed threshold \(\beta\).

In this nonlinear setting, an exact covering description is not apparent from the analytic form of~\eqref{eq:sigmoid}. To the author's knowledge, such a covering reformulation has not been proposed for constraints of this form. 
The proposed framework derives the reformulation by shifting attention from the analytic form of \(f\) to the order structure of its feasible binary supports.
Since \(f\) is monotone in \(x\), the feasible set \(\Omega=\{S\subseteq[n]:f(x_S;\hat w,\hat\alpha)\ge\beta\}\) is upper-closed. Hence the covering formulation~\eqref{eq:pformub} applies, indexed by the dual structure set \(\ec\nc\Omega\). Then, a direct computation of \(\ec\nc\Omega\) gives
\[
  \ec\nc\Omega = \{T \subseteq \Delta:=[n] \mid f(x_{\Delta\setminus T};\hat w,\hat\alpha)<\beta\}.
\]
Moreover, since membership can be verified efficiently by evaluating the function, Corollary~\ref{coro:sep} implies that, for any binary $x$, the associated inequality can be separated in polynomial time.
\hfill $\triangle$
\end{example}
\end{revision}

\subsection{\rev{Generalizing Polyhedral Analysis}}
\label{sec:polyhedral}
Another unifying perspective offered by Theorem~\ref{thm:reform} is that if either the solution space $\Omega$ or the objective function $f$ is monotone, the binary optimization problem $\min_{x \in \mathcal{X}_\Omega} f(x)$ can be reformulated exactly as a set covering formulation \eqref{eq:pformub} or an elimination formulation \eqref{eq:pformlb}. 
\rev{
This reformulation bridges general BIPs with the classical set covering literature, allowing polyhedral tools from set covering to be transferred to broader BIP settings. We illustrate this transfer through facet analysis in this subsection and inequality strengthening in the next.
}
To unify the facet analysis of covering and elimination inequalities, we first introduce the following facet-preserving \rev{lemma}.

\begin{restatable}{lemma}{faceti}
  \label{lem:faceti}
Given a binary solution space $\mathcal{X} := \{x \in \{0,1\}^n \mid g(x) \leq 0\}$ and an index set $I \subseteq [n]$, define the flipping map $\theta_I : \{0,1\}^n \to \{0,1\}^n$ by
\[
(\theta_I(x))_i := \begin{cases}
1 - x_i & \text{if } i \in I, \\
x_i & \text{otherwise}.
\end{cases}
\]
Let $\mathcal{X}_{\theta_I} := \{x \in \{0,1\}^n \mid g(\theta_I(x)) \leq 0\}$. Then, $\theta_I$ induces a bijective affine map between $\conv(\mathcal{X})$ and $\conv(\mathcal{X}_{\theta_I})$. In particular, an inequality $\langle a, x \rangle \geq b$ is valid (or facet-defining) for $\conv(\mathcal{X})$ if and only if $\langle a, \theta_I(x) \rangle \geq b$ is valid (or facet-defining) for $\conv(\mathcal{X}_{\theta_I})$.
\end{restatable}

Since the elimination problem~\eqref{eq:pformlb} can be reformulated as a covering problem~\eqref{eq:pformub} via the flipping map \(\theta_{[n]}\), its facet analysis can be unified with that of the covering problem. This equivalence allows general techniques from the set covering literature, including facet-defining conditions~\citep{balas1972set,balas1989setii,sanchez1998set} and lifting-based constraint strengthening~\citep{wei2022integer,nobili1989facets,wei2025supervalid}, to be extended to a broader class of binary integer programs. 
While a comprehensive treatment is beyond the scope of this paper, 
we illustrate this through a facet-defining criterion for general BIPs adapted from the existing \rev{set covering} literature~\citep{balas1972set,wei2022integer}, which requires the following definition.

\begin{definition}[Quasi-Feasibility]
  A structure $T \in \nc\Omega$ is termed \emph{quasi-feasible} if, for every element $a \in T$, there exists an element $a' \in \Delta \setminus T$ such that $T \setminus \{a\} \cup \{a'\} \in \Omega$.
\end{definition}

This concept identifies infeasible sets that are locally close to feasibility and is used to state facet-defining conditions for covering-type constraints. \rev{We emphasize that the value is not in rederiving these classical conditions, which are well known in the set covering literature, but in using the proposed framework as a bridge to transfer them to general BIPs with monotone structures.}


\begin{restatable}{proposition}{facetthm}
  \label{prop:facet}
  Suppose $\Omega$ is an upper system and $|T| \geq 2$ for every $T \in m\ec{\nc{(\Omega)}}$, the covering inequality \eqref{eq:pformub01} is facet-defining if and only if $\Delta \setminus T\in M\nc(\Omega)$ is quasi-feasible.
\end{restatable}

\begin{revision}
By Lemma~\ref{lem:faceti}, the elimination constraints in~\eqref{eq:pformlb01} can be viewed as covering constraints after applying the flipping map. Specifically, setting \(x'_i=1-x_i\) gives
\[
\sum_{i \in T} x_i \leq |T|-1, \ \forall T \in m\nc\Omega \quad \Longleftrightarrow \quad \sum_{i \in T} x'_i \geq 1, \ \forall T\in m\nc\Omega = m\ec\nc(\ec\Omega).
\]
Thus, the preceding covering-side results can be applied to a lower system \(\Omega\) by considering its flipped upper system \(\ec\Omega\). This yields the following lower-system analogue.

\begin{restatable}{proposition}{facetthmlower}
  \label{prop:facetlower}
  Suppose $\Omega$ is a lower system and $|T| \geq 2$ for every $T \in m{\nc{(\Omega)}}$, the elimination inequality \eqref{eq:pformlb01} is facet-defining if and only if $\Delta \setminus T\in M\nc(\ec\Omega)$ is quasi-feasible with respect to $\ec\Omega$. That is, for every $a \in \Delta \setminus T$ there is some $a' \in T$ such that $T \setminus \{a'\} \cup \{a\} \in \Omega$.
\end{restatable}



These examples illustrate that the framework provides a systematic way to translate covering-side polyhedral results into facet conditions for BIP-derived structures.
We demonstrate this through two applications: (i) 
a new facet characterization for the bilinear formulation introduced in Example~\ref{eg:bilinear};
(ii) a simplified facet proof for the bipartite induced subgraph problem.

\begin{example}[Linearization of Bilinear Terms (Continued)]
  \label{eg:bilinearcont}
  Example~\ref{eg:bilinear} derives an exact reformulation for the structured bilinear program with the following covering constraints
    $$\sum_{i \in I} x_i + \sum_{j \in J} y_j \geq 1, \quad \forall (I,J) \in m\ec\nc(\Omega),$$
    where $\Omega:=\{(I, J) \mid \sum_{(i,j) \in I\times J} R_{ij} \geq \alpha\}$. Then, 
    Proposition~\ref{prop:facet} gives a one-line facet proof for the basic formulation without additional constraints.
Suppose every submatrix obtained by removing only one row or one column from \(R\) has total weight at least \(\alpha\). Then the above constraint is facet-defining if and only if: (i) \((I,J)\) is minimal, in the sense that adding back any row from \(I\) or any column from \(J\) makes the remaining submatrix feasible; and (ii) for every row or column in the submatrix \(R_{[n]\setminus I,\,[m]\setminus J}\), there exists a row in \(I\) or a column in \(J\) whose swap yields a feasible submatrix.
\hfill $\triangle$
\end{example}

\begin{example}[Bipartite Induced Subgraph Problem]
Given a graph \(G=(V,E)\) with nonnegative vertex weights, the bipartite induced subgraph problem (BISP) seeks a maximum-weight vertex set \(T\subseteq V\) such that \(G[T]\) is bipartite. This problem has been studied under several names in the literature~\citep{fouilhoux2006polyhedral,barahona1989facets}. The associated structure set
\(
\Omega=\{T\subseteq V \mid G[T]\text{ is bipartite}\}
\)
is lower-closed and therefore admits an exact elimination reformulation indexed by \(m\nc\Omega\). Since \(\nc\Omega\) consists of vertex sets whose induced subgraphs contain odd cycles, \(m\nc\Omega\) consists precisely of vertex sets \(T\) for which \(G[T]\) is an odd cycle. Hence, the problem can be written as
\[
\max_{x\in\{0,1\}^{|V|}}
\left\{
\iprod{w,x}
~\middle|~
\sum_{i\in T}x_i\le |T|-1,\ \forall T\in m\nc\Omega
\right\}.
\]
Since every such \(T\) satisfies \(|T|\ge 2\), Proposition~\ref{prop:facetlower} implies that the inequality associated with \(T\) is facet-defining if and only if, for every \(v\in V\setminus T\), there exists \(v'\in T\) such that \(G[T\setminus\{v'\}\cup\{v\}]\) is bipartite. Equivalently, after adding \(v\) to the odd cycle \(G[T]\), one must be able to delete a single vertex of \(T\) to destroy all odd cycles. This occurs exactly when \(G[T\cup\{v\}]\) contains at most two odd cycles, recovering the condition derived in~\citet[Theorem 2.7]{barahona1989facets} through the general quasi-feasibility criterion.
\hfill \(\triangle\)
\end{example}
\end{revision}

\begin{revision}
\subsection{Transferring Lifted Inequalities} 
\label{sec:ineq}
Besides supporting facet analysis for basic covering inequalities, the covering reformulation in Theorem~\ref{thm:reform} also allows strengthened inequalities from the set covering literature to be transferred to general BIPs. This subsection illustrates this through the \emph{vitality inequalities} introduced in~\citet{wei2022integer}, which generate valid inequalities from any set \(S\subseteq\Delta\) that contains a structure \(T\in\ec\nc\Omega\).
The following definitions are adapted from~\citet{wei2022integer}. The lower-system version follows from the upper-system case by applying the transformation \(\Omega \mapsto \ec\Omega\).

\begin{definition}[Vitality Inequalities (Upper Systems)]
  \label{defi:vitaupper}
  Given $S \subseteq \Delta$, its vital components are defined as $\mathcal V_S = \{S' \subseteq S \mid (\Delta \setminus S) \cup S' \in \Omega\}$. Let $m(\mathcal V_S)$ denote the minimal elements in $\mathcal V_S$, a disjoint decomposition is a partition $\{\mathcal V_i\}_{i \in I}$ of $m(\mathcal V_S)$ such that $(\bigcup \mathcal V_i) \cap (\bigcup \mathcal V_j) = \emptyset$ for $i\neq j\in I$. 
  For each part $\mathcal V_i$, let $\gamma_i:= \min_{S' \in \mathcal V_i}|S'|$. Then the associated vitality inequality is 
  \[\sum_{i \in I}\left(\prod_{j\neq i} \gamma_j\right) \sum_{a \in \bigcup\mathcal V_i} x_a \geq \prod_{i \in I} \gamma_i.\]
\end{definition}

\begin{definition}[Vitality Inequalities (Lower Systems)]
  \label{defi:vitalower}
  Given \(S \subseteq \Delta\), its vital components are defined as
  \(
  \mathcal V_S = \{S' \subseteq S \mid S \setminus S' \in \Omega\}.
  \)
  Let $m(\mathcal V_S)$ denote the minimal elements in $\mathcal V_S$, a disjoint decomposition is a partition $\{\mathcal V_i\}_{i \in I}$ of $m(\mathcal V_S)$ such that $(\bigcup \mathcal V_i) \cap (\bigcup \mathcal V_j) = \emptyset$ for $i\neq j\in I$. 
  For each part \(\mathcal V_i\), let
  \(
  \gamma_i := \min_{S'\in \mathcal V_i}|S'|.
  \)
  Then the associated vitality inequality is
  \[
  \sum_{i \in I}
  \left(\prod_{j\neq i} \gamma_j\right)
  \sum_{a \in \bigcup \mathcal V_i} x_a
  \leq
  \sum_{i \in I}
  \left(\prod_{j \neq i} \gamma_j\right)
  \left|\bigcup \mathcal V_i\right|
  -
  \prod_{i \in I} \gamma_i.
  \]
\end{definition}

In both definitions, \(\mathcal V_S\) represents the family of local modifications of \(S\) that restore feasibility.
For general $S$, enumerating $\mathcal{V}_S$ may be intractable (e.g., $S = \Delta$), but when $S$ has sufficient structure for $\mathcal{V}_S$ to be explicitly characterized, the associated vitality inequality can strengthen the formulation.
Many classical strengthened inequalities in graph optimization arise this way. 
For example, consider the lower-closure of Hamiltonian cycles (HCs) in a complete graph, $\Omega = \{T \subseteq \Delta \mid T \text{ forms a linear forest or is a HC}\}$, and let $S$ be a clique of size $k < |V|$. Then $\mathcal{V}_S = \{S' \subseteq S \mid S \setminus S' \text{ forms a linear forest}\}$, and without any disjoint decomposition, $\gamma = \binom{k}{2} - (k-1)$ gives the minimum number of edges to remove to ensure acyclicity within $S$. The vitality inequality from Definition~\ref{defi:vitalower} then recovers exactly the subtour elimination inequalities. 
The following two examples show how this construction can be transferred from the covering literature to general BIPs.

\begin{example}[Linearization of Bilinear Terms (Continued)]
  \label{eg:bilift}
  Consider a minimal \((I,J)\in m\ec\nc\Omega\), so that
  \(
  \sigma(\bar I,\bar J)<\alpha.
  \)
  Define
  \(
  g:=\alpha-\sigma(\bar I,\bar J)>0
  \) as the feasibility gap.
  Now choose a row \(i'\notin I\) and form the larger set
  \(
  S=(I\cup\{i'\},J).
  \)
  The new gap becomes \(g'=g+\sigma(\{i'\},\bar J)\). By definition, $\mathcal V_S$ consists of all index sets $(I',J')$ with $I' \subseteq I \cup \{i'\}$ and $J' \subseteq J$ such that $\sigma( \overline{I \cup \{i'\} \setminus I'}, \overline{J \setminus J'}) \geq \alpha$. To construct a disjoint decomposition of $m(\mathcal V_S)$, call a row \(i\in I\) (or a column $j \in J$) \emph{strong} if
  \(
  \sigma(\{i\},\bar J)\ge g'
  \) (or $\sigma(\bar I\setminus\{i'\},\{j\})\ge g'$).
  Let \(I_0\subseteq I\) and \(J_0\subseteq J\) denote the strong rows and columns. The strong elements form the first part with $\gamma_1 = 1$, and the remaining minimal elements in $\mathcal V_S$ form the second part with $\gamma_2 = 2$.
  Hence, by Definition~\ref{defi:vitaupper}, we obtain the vitality inequality
    \begin{equation}
      \label{eq:vat}
  2\left(\sum_{i\in I_0}x_i+\sum_{j\in J_0}y_j\right)
  +x_{i'}
  +\sum_{i\in I\setminus I_0}x_i
  +\sum_{j\in J\setminus J_0}y_j
  \ge 2.
    \end{equation}
  If no strong rows or columns are identified, the construction reduces to the single-part vitality inequality with right-hand side \(2\).
\hfill \(\triangle\)
\end{example}

\citet{barahona1989facets} identified basic odd-cycle elimination cuts and their lifted clique variants for the maximum bipartite induced subgraph problem, and \citet{fouilhoux2006polyhedral} later introduced valid inequalities based on two special graph structures termed wheels and path-cycles. The following example derives another class of valid inequalities from the vital components in Definition~\ref{defi:vitalower}, which, to the author's knowledge, is new for this problem.

\begin{example}[Bipartite Induced Subgraph Problem (Continued)]
Consider any \(S \subseteq V\) such that \(G[S]\) forms \(K\) internally vertex-disjoint \(s\)-\(t\) paths \(\mathcal P=\{P_k\}_{k \in [K]}\) with two shared endpoints \(s,t \in V\) (see Figure~\ref{fig:st-parity-paths}). We partition \(\mathcal P\) into \(\gamma_{\mathrm{odd}}\) odd paths and \(\gamma_{\mathrm{even}}\) even paths, depending on the parity of the number of edges in each \(s\)-\(t\) path. Then, every odd cycle in \(G[S]\) is induced by the combination of one odd path with one even path. Three types of minimal \(S' \subseteq S\) can make \(G[S \setminus S']\) odd-cycle-free. First, \(\mathcal V_1 = \{\{s\}, \{t\}\}\), with \(\gamma_1=1\). Second, one can select one internal vertex from each odd path, with solution size \(\gamma_{\mathrm{odd}}\). Third, one can select one internal vertex from each even path, with solution size \(\gamma_{\mathrm{even}}\). Clearly, these three groups have no intersection in vertices, forming a disjoint decomposition. Then, the associated lower-system vitality inequality is
\[
\sum_{i \in \{s,t\}} \gamma_{\mathrm{odd}}\gamma_{\mathrm{even}}x_i
+ \sum_{i \in V_{\mathrm{odd}}} \gamma_{\mathrm{even}} x_i
+ \sum_{i \in V_{\mathrm{even}}} \gamma_{\mathrm{odd}} x_i
\leq
\gamma_{\mathrm{even}}|V_{\mathrm{odd}}|
+ \gamma_{\mathrm{odd}}|V_{\mathrm{even}}|
+ \gamma_{\mathrm{odd}}\gamma_{\mathrm{even}},
\]
where \(V_{\mathrm{odd}}\) and \(V_{\mathrm{even}}\) denote the internal vertices of the odd and even paths, respectively. 
When both \(\gamma_{\mathrm{odd}}\) and \(\gamma_{\mathrm{even}}\) are strictly greater than \(1\), the above inequality can be further strengthened by applying the Chv\'{a}tal-Gomory rounding procedure.
\hfill \(\triangle\)
\end{example}

\begin{figure}[tb]
\centering
\begin{tikzpicture}[
    scale=1,
    every node/.style={font=\small},
    endpoint/.style={circle, draw, fill=black, inner sep=1.8pt},
    vnode/.style={circle, draw, fill=black, inner sep=1.4pt},
    lab/.style={font=\footnotesize},
    thick
]
\pgfmathsetmacro{\gap}{0.9}

\node[endpoint, label=left:$s$] (s) at (0,0) {};
\node[endpoint, label=right:$t$] (t) at (10,0) {};
\node[vnode] (a1) at (2,2*\gap) {};
\node[vnode] (a2) at (4,2*\gap) {};
\node[vnode] (a3) at (6,2*\gap) {};
\node[vnode] (a4) at (8,2*\gap) {};
\draw (s) -- (a1) -- (a2) -- (a3) -- (a4) -- (t);
\node[lab] at (5,2*\gap+0.3) {$P_1$ (odd)};
\node[vnode] (b1) at (2.8,\gap) {};
\node[vnode] (b2) at (7.2,\gap) {};
\draw (s) -- (b1) -- (b2) -- (t);
\node[lab] at (5,\gap+0.3) {$P_2$ (odd)};
\node[vnode] (c1) at (5,0) {};
\draw (s) -- (c1) -- (t);
\node[lab] at (5,0.3) {$P_3$ (even)};
\node[vnode] (d1) at (2.5,-\gap) {};
\node[vnode] (d2) at (5.0,-\gap) {};
\node[vnode] (d3) at (7.5,-\gap) {};
\draw (s) -- (d1) -- (d2) -- (d3) -- (t);
\node[lab] at (5,-\gap+0.3) {$P_4$ (even)};
\node[vnode] (e1) at (2.2,-2*\gap) {};
\node[vnode] (e2) at (5.0,-2*\gap) {};
\node[vnode] (e3) at (7.8,-2*\gap) {};
\draw (s) -- (e1) -- (e2) -- (e3) -- (t);
\node[lab] at (5,-2*\gap+0.3) {$P_5$ (even)};
\end{tikzpicture}
\caption{A collection of five internally vertex-disjoint \(s\)-\(t\) paths, with two odd paths and three even paths.}
\label{fig:st-parity-paths}
\end{figure}

\end{revision}

\section{Reformulations via Latent Monotone Systems}
\label{sec:bipreform}
The exactness of reformulations \eqref{eq:pformub}--\eqref{eq:pformlb} requires monotonicity of either the objective or feasible region. 
\rev{
This section extends their scope to BIPs whose monotone structure is latent rather than explicit, through two reformulation strategies.
Although these reformulations follow naturally once the transformations are presented, they remain underexplored because covering and elimination inequalities are typically treated as structure-specific cutting planes, 
with their connection to monotone set systems left implicit.
This section explores latent monotone structures to extend the usage of covering and elimination inequalities.
}

\subsection{Bimonotone Cuts}
\label{sec:bimono}
The first source of latent monotone structure comes from functions that are monotone with opposite signs on a bipartition of the variables. We need the following \rev{definition}, where a bipartition is said to be trivial if one of the two parts is empty.

\begin{revision}
\begin{definition}[Bimonotone Function]
  A function \(g:\{0,1\}^\Delta\to\mathbb R\) is called bimonotone if there exists a possibly trivial bipartition \((I,J)\) of \(\Delta\) such that, for every fixed \(x_J\), \(g(\cdot,x_J)\) is increasing in \(x_I\), and for every fixed \(x_I\), \(g(x_I,\cdot)\) is decreasing in \(x_J\).
\end{definition}

The following proposition identifies several common classes of bimonotone functions.
\end{revision}

\begin{restatable}{proposition}{func}
  \label{prop:func}
  The following functions are bimonotone:
  \begin{itemize}
    \item Linear (modular) functions.
    \item Bilinear functions $\iprod{x, Ry}$ for some block-diagonal matrix $R \in \mathbb R^{I \times J}$ where there exist (possibly trivial) index bipartitions $I= I_1 + I_2$ and $J = J_1 + J_2$ such that the diagonal blocks satisfy $R_{I_1J_1} \geq 0$ and $R_{I_2J_2} \leq 0$ (entrywise).
    \item Submodular functions $f$ where for every $i \in \Delta$, either 
      $$f(\{i\})-f(\emptyset) \leq 0 \text{ or } f(\Delta) - f(\Delta \setminus \{i\}) \geq 0.$$
    \item The supermodular counterpart.
  \end{itemize}
\end{restatable}

Once such a bimonotone objective function is identified in a given BIP, we obtain an exact reformulation method and facet-defining conditions according to the following corollary.

\begin{revision}
\begin{restatable}{corollary}{bimono}
  \label{coro:bimono}
  Given a BIP \(\min_{x=(x_I,x_J)\in \mathcal X_{\Omega}} f(x_I,x_J)\) where \(f\) is bimonotone under the index partition \((I,J)\). Let \(\theta(x_I,x_J)=(x_I,1-x_J)\) be the flipping map on \(J\), and define the transformed system
  \(
    \Omega_\theta
    := \{(T_I,J\setminus T_J) \mid (T_I,T_J)\in\Omega\}.
  \)
  Then the problem admits the following reformulation:
  \begin{equation}
    \label{eq:bimono}
    \begin{aligned}
      \min\limits_{x\in\{0,1\}^n} &~ f(x) \\
      \text{s.t.} 
      &~ \sum_{i\in T_I} x_i
      + \sum_{j \in T_J}(1-x_j)
      \ge 1,
      \quad \forall (T_I, T_J)\in m\ec{\nc{\uparrow\Omega_\theta}} .
    \end{aligned}
  \end{equation}
Moreover, suppose $|T_I \cup T_J| \geq 2$ for every $(T_I, T_J) \in m\ec{\nc{\uparrow\Omega_\theta}}$, the above inequality indexed by $(T_I, T_J)$ is facet-defining if and only if $(I \setminus T_I, J \setminus T_J)$ is quasi-feasible with respect to $\uparrow\Omega_\theta$.
\end{restatable}

\end{revision}




These results can help analyze the structures tied to \eqref{eq:bimono} and separate valid inequalities or even facets to strengthen the associated BIP. We demonstrate this with the following example.

\begin{revision}
\begin{example}[Min-Cut with Signed Weights]
  \label{eg:signmincut}
Consider a graph \(G=(V,E)\) with signed edge weights \(\{w_e\}_{e\in E}\). The signed min-cut problem seeks a minimum-weight edge cut. With mixed signs, the associated Laplacian is generally indefinite, so the convex QP formulation available in the nonnegative case no longer applies.
By Corollary~\ref{coro:bimono}, we partition the edge set as \(E=(I:=E_+,J:=E_-)\), where \(E_+\) and \(E_-\) contain the nonnegative- and negative-weight edges. The objective function is bimonotone with respect to this bipartition, and formulation \eqref{eq:bimono} yields an exact linear representation of the signed min-cut problem.
In this setting, \(\uparrow\Omega_\theta\) consists of edge subsets \((T_I,J\setminus T_J)\) such that \(T_I\supseteq C_I\) and \(T_J\subseteq C_J\) for some edge cut \((C_I,C_J)\). Equivalently, edges in \(I\setminus T_I\) must remain within one side of the cut, while edges in \(T_J\) must cross the cut. Hence, membership of \((T_I,J\setminus T_J)\) in \(\uparrow\Omega_\theta\) can be tested by contracting each connected component of the subgraph \((V,I\setminus T_I)\), retaining only the projected edges in \(T_J\), and checking whether the resulting graph is bipartite.
Thus, the reformulation yields a BIP with linear inequalities whose integer point separation is supported by a polynomial-time membership oracle. By Corollary~\ref{coro:sep}, this oracle enables efficient integer point separation for~\eqref{eq:bimono}, and hence supports a cut-generation implementation.
\hfill $\triangle$
\end{example}
\end{revision}

\subsection{Interval System Decomposition}

In addition to bimonotone functions, latent monotone structure can also be revealed through set systems that exhibit the interval property, as characterized below.

\begin{definition}[Interval System]
A set system $\Omega$ is said to have the interval property, or to be an interval system, if for all $T, T' \in \Omega$ and any $T''$ with $T \subseteq T'' \subseteq T'$, it follows that $T'' \in \Omega$.
\end{definition}

\begin{restatable}{proposition}{interval}
  \label{lem:interval}
  For any set system $\Omega$, $\uparrow\Omega \cap \downarrow\Omega$ is the smallest interval system that contains $\Omega$. Moreover, $\Omega$ is interval if and only if $\Omega= \uparrow\Omega \cap \downarrow\Omega$.
\end{restatable}

Therefore, every identified interval system $\Omega$ can be represented as follows
  $$\mathcal X_\Omega = \left\{x \in \{0,1\}^n \;\middle|\; \sum_{i \in T}x_i \geq 1, \forall T \in m\ec{\nc{\uparrow(\Omega)}}, \sum_{i \in T}x_i \leq |T| - 1, \forall T \in m{\nc{\downarrow(\Omega)}} \right\}.$$
  \begin{revision}
  Accordingly, for every interval system decomposition $\Omega = \bigcup_{k \in K} \Omega_k$ where each $\Omega_k$ is an interval system, the following reformulation is exact.
  \begin{equation}
    \label{eq:decomp}
  \begin{aligned} 
    \min\limits_{x \in \{0,1\}^n, z \in \{0,1\}^{|K|}} &~ f(x)\\
    \text{s.t.} & ~\sum_{i \in T} x_i \geq z_k, ~~~~~~~~\forall k \in K, \forall T \in m\ec{\nc{\uparrow(\Omega_k)}},\\
                & ~\sum_{i \in T} x_i \leq |T| - z_k, ~\forall k \in K, \forall T \in m\nc{\downarrow(\Omega_k)}, \\
                & ~\sum_{k \in K} z_k = 1.
\end{aligned}
\end{equation}
  \end{revision}
  

  Although such a decomposition always exists (e.g., trivial singleton intervals $\Omega=\bigcup_{T \in \Omega} \uparrow\{T\} \cap \downarrow\{T\}$), its practical value arises when the number of components is small or structured, a situation that naturally occurs in disjunctive-type problems \citep{balas1979disjunctive}.



\begin{revision}
\begin{example}[Piecewise Monotone Objectives]
  \label{eg:piecewise}
In network design, one may wish to minimize a surrogate objective $f(x)$ learned from data \citep{ma2015large,elfar2018machine}. For example, tree-based models \citep{ke2017lightgbm,loh2011classification} yield a piecewise-affine predictor with each leaf defining a region $\mathcal X_k$ and the objective $f(x) = \iprod{c^k, x} + c_0^k$.
Suppose each $c^k$ is positive (index set $K_+$) or negative ($K_-$). Then $\min_{x \in \mathcal X_\Omega} f(x)$ admits the following exact reformulation
  $$
   \begin{aligned}
     \min\limits_{x \in \{0,1\}^n, z \in \{0,1\}^{|K|}} &~ \eta \\
     \text{s.t.} & ~\eta \geq \iprod{c^k, x} + c_0^k - M(1-z_k), ~~~~\quad\quad\quad\quad\quad\forall k \in K,\\
                 & ~\sum_{i \in T} x_i \geq z_k, ~~~~~~~~\forall k \in K_+, \forall T \in m\ec{\nc{\uparrow{(\Omega_{\mathcal X_k})}}},\\
                 & ~\sum_{i \in T} x_i \leq |T| - z_k, ~\forall k \in K_-, \forall T \in m\nc{\downarrow{(\Omega_{\mathcal X_k})}},\\
                 & ~\sum_{k \in K} z_k = 1,
  	\end{aligned}
  $$
  where $M>0$ is sufficiently large. The binary variable $z_k$ selects the active region. For $k \in K_+$ the covering inequalities enforce $x \in \uparrow\Omega_{\mathcal X_k}$; for $k \in K_-$ the elimination inequalities apply.
\hfill $\triangle$
\end{example}
\end{revision}

\rev{
These results show that identifying latent monotone systems broadens the scope of monotone reformulations by connecting them to bimonotone, disjunctive, and piecewise structures.
}

\section{Computational Experiments}
\label{sec:case}
\begin{revision}
This section illustrates the proposed framework through two case studies: (i) a bilinear supplier-region selection problem, where bimonotone cuts and their lifted valid inequalities are tested both as a standalone exact approach and as root-node strengthening cuts; and (ii) a distributionally robust network site selection problem, where interval set systems lead to flexible cut-based reformulations.
All experiments were conducted on a 2023 MacBook Pro with an M2 Max chip, 12 CPU cores, and 64 GB of memory. The implementations use Python 3.10 and Gurobi 13.0.2 as the solver.
\end{revision}

\begin{revision}

\subsection{Case Study 1: Bilinear Supplier--Region Selection}
\label{sec:case1}

Bilinear binary programs with nonnegative interaction matrices arise in a range of combinatorial settings, including quadratic assignment~\citep{koopmans1957assignment}, hub location~\citep{o1987quadratic}, and quadratic knapsack~\citep{gallo2009quadratic}, where the term $\iprod{x, Ry}$ with $R \geq 0$ captures synergy, coverage, or compatibility from jointly selecting two classes of decisions.
To demonstrate the proposed framework in a bimonotone setting, we extend the monotone supplier--region selection problem of Example~\ref{eg:bilinear} by allowing both positive and adverse interaction blocks. Let $(I, \bar{I})$ and $(J, \bar{J})$ partition the suppliers and regions, respectively, where interactions in $I \times J$ generate positive partnership value and those in $\bar{I} \times \bar{J}$ incur adverse effects such as incompatibility. This yields
\begin{equation}
  \label{eq:bibilinear}
\begin{aligned}
\min_{x\in\{0,1\}^n,\;y\in\{0,1\}^m}
&\quad \iprod{c_1,x}+\iprod{c_2,y} \\
\text{s.t.}
&\quad
\iprod{x_I, R_{IJ}\,y_J}
-
\iprod{x_{\bar I}, R_{\bar I\bar J}\,y_{\bar J}}
\geq \alpha,
\end{aligned}
\end{equation}
where $R_{IJ}, R_{\bar{I}\bar{J}} \geq 0$ and entries in $c_1, c_2$ may take either sign.

We conduct two experiments to evaluate the computational impact of the proposed bimonotone cuts and their lifted variants. Let \(n=m\) denote the dimension of $x$ and $y$. For each configuration \((n,\gamma)\), we randomly generate three instances, where \(c_1,c_2\in\mathbb R^{n}\) are drawn uniformly from \((-10,10)\), and \(R_{IJ},R_{\bar I\bar J}\ge 0\) are \(n/2\times n/2\) matrices with entries drawn from a mixture distribution: \(U(1,5)\) with probability \(0.8\) and \(U(20,50)\) with probability \(0.2\). We set \(\alpha:=\gamma\sum_{i\in I,j\in J}R_{ij}\), ensuring feasibility for the tested \(\gamma\) values. Results are averaged over three instances. We report runtime and metrics related to separation and root-node cut generation, whose core routines are compiled with Numba and parallelized for efficiency.

\subsubsection{Exact Approach via Bimonotone and Lifted Cuts}
This experiment compares four solution approaches.
\textsc{Vanilla} solves the bilinear program~\eqref{eq:bibilinear} directly using Gurobi's native bilinear solver.
\textsc{Lin} solves the following McCormick linearization of~\eqref{eq:bibilinear} with auxiliary variables $z_{ij}=x_i y_j$ indexed by $I\times J$ and $\bar I\times \bar J$:
\[
\begin{aligned}
\min_{x\in\{0,1\}^n,\;y\in\{0,1\}^m,\;z\in[0,1]}
&\quad \iprod{c_1,x}+\iprod{c_2,y} \\
\text{s.t.}
&\quad
\iprod{R_{IJ},z_{IJ}}
-
\iprod{R_{\bar I\bar J},z_{\bar I\bar J}}
\geq \alpha,\\
&\quad z_{ij} \leq x_i,\quad z_{ij} \leq y_j,
\quad \forall i\in I,\ j\in J,\\
&\quad z_{ij} \geq x_i+y_j-1,
\quad \forall i\in \bar I,\ j\in \bar J.
\end{aligned}
\]
\textsc{Cut} implements the pure bimonotone cut reformulation~\eqref{eq:bimono}, whereas \textsc{LiftCut} generates the vitality inequality~\eqref{eq:vat} when available and falls back to the standard bimonotone cut otherwise.
We vary $n$ over $\{10,15,\ldots,35\}$ and $\gamma$ over $\{0.7,0.75,\cdots,0.95\}$, with three random instances per configuration, yielding $108$ instances in total. We impose a time limit of $300$ seconds for each run.

\vbadness=10001
\begin{table}[!tbp]
  \centering
\rev{
 \scriptsize
\setlength{\tabcolsep}{3pt}
\begin{tabular}{lrrrr@{\hskip 10pt}rr@{\hskip 10pt}rr}
\toprule
Config & \multicolumn{4}{c}{Solve Time (s)}
& \multicolumn{2}{c}{Sep. Time (s)}
& \multicolumn{2}{c}{Num. Cuts}\\
\cmidrule(lr){1-1}\cmidrule(lr){2-5}\cmidrule(lr){6-7}\cmidrule(lr){8-9}
$(n,\gamma)$ & Vanilla & Lin & Cut & LiftCut & Cut & LiftCut & Cut & LiftCut \\
\midrule
$(10, 0.7)$ & $0.002$ & $\mathbf{0.001}$ & $0.001$ & $0.002$ & $\mathbf{0.000}$ & $0.000$ & $16.7$ & $\mathbf{16.3}$ \\
$(10, 0.75)$ & $0.015$ & $0.006$ & $0.014$ & $\mathbf{0.002}$ & $\mathbf{0.000}$ & $0.000$ & $14.3$ & $\mathbf{13.7}$ \\
$(10, 0.8)$ & $0.001$ & $0.001$ & $\mathbf{0.001}$ & $0.001$ & $\mathbf{0.000}$ & $0.000$ & $13.3$ & $\mathbf{13.0}$ \\
$(10, 0.85)$ & $0.001$ & $\mathbf{0.000}$ & $0.001$ & $0.001$ & $\mathbf{0.000}$ & $0.000$ & $13.0$ & $\mathbf{12.7}$ \\
$(10, 0.9)$ & $0.001$ & $\mathbf{0.000}$ & $0.000$ & $0.001$ & $\mathbf{0.000}$ & $0.000$ & $\mathbf{11.7}$ & $\mathbf{11.7}$ \\
$(10, 0.95)$ & $0.001$ & $\mathbf{0.001}$ & $0.001$ & $0.001$ & $\mathbf{0.000}$ & $0.000$ & $12.0$ & $\mathbf{11.7}$ \\
\cmidrule(lr){1-9}
$(15, 0.7)$ & $0.003$ & $\mathbf{0.002}$ & $0.008$ & $0.009$ & $\mathbf{0.000}$ & $0.001$ & $\mathbf{62.7}$ & $81.0$ \\
$(15, 0.75)$ & $0.003$ & $\mathbf{0.002}$ & $0.016$ & $0.019$ & $\mathbf{0.000}$ & $0.000$ & $\mathbf{51.0}$ & $51.7$ \\
$(15, 0.8)$ & $0.003$ & $\mathbf{0.002}$ & $0.011$ & $0.009$ & $\mathbf{0.000}$ & $0.000$ & $46.0$ & $\mathbf{45.0}$ \\
$(15, 0.85)$ & $0.003$ & $\mathbf{0.002}$ & $0.004$ & $0.007$ & $\mathbf{0.000}$ & $0.000$ & $\mathbf{29.3}$ & $33.7$ \\
$(15, 0.9)$ & $0.002$ & $\mathbf{0.001}$ & $0.004$ & $0.003$ & $\mathbf{0.000}$ & $0.000$ & $26.7$ & $\mathbf{21.7}$ \\
$(15, 0.95)$ & $0.002$ & $0.001$ & $\mathbf{0.001}$ & $0.001$ & $\mathbf{0.000}$ & $0.000$ & $18.0$ & $\mathbf{17.3}$ \\
\cmidrule(lr){1-9}
$(20, 0.7)$ & $0.003$ & $\mathbf{0.002}$ & $0.047$ & $0.046$ & $0.001$ & $\mathbf{0.001}$ & $297.0$ & $\mathbf{232.7}$ \\
$(20, 0.75)$ & $0.004$ & $\mathbf{0.003}$ & $0.031$ & $0.041$ & $0.001$ & $\mathbf{0.001}$ & $309.7$ & $\mathbf{217.0}$ \\
$(20, 0.8)$ & $0.004$ & $\mathbf{0.004}$ & $0.032$ & $0.024$ & $\mathbf{0.001}$ & $0.001$ & $186.3$ & $\mathbf{184.3}$ \\
$(20, 0.85)$ & $0.004$ & $\mathbf{0.003}$ & $0.041$ & $0.023$ & $\mathbf{0.000}$ & $0.001$ & $103.7$ & $\mathbf{101.3}$ \\
$(20, 0.9)$ & $0.004$ & $\mathbf{0.002}$ & $0.022$ & $0.032$ & $\mathbf{0.000}$ & $0.000$ & $90.0$ & $\mathbf{81.3}$ \\
$(20, 0.95)$ & $0.002$ & $\mathbf{0.002}$ & $0.005$ & $0.005$ & $\mathbf{0.000}$ & $0.000$ & $\mathbf{36.3}$ & $\mathbf{36.3}$ \\
\cmidrule(lr){1-9}
$(25, 0.7)$ & $0.006$ & $\mathbf{0.005}$ & $0.151$ & $0.141$ & $0.006$ & $\mathbf{0.005}$ & $1425.3$ & $\mathbf{907.0}$ \\
$(25, 0.75)$ & $0.006$ & $\mathbf{0.004}$ & $0.107$ & $0.107$ & $0.005$ & $\mathbf{0.004}$ & $1136.7$ & $\mathbf{726.0}$ \\
$(25, 0.8)$ & $0.007$ & $\mathbf{0.005}$ & $0.091$ & $0.097$ & $\mathbf{0.003}$ & $0.003$ & $643.3$ & $\mathbf{519.7}$ \\
$(25, 0.85)$ & $0.006$ & $\mathbf{0.005}$ & $0.046$ & $0.052$ & $0.002$ & $\mathbf{0.001}$ & $428.0$ & $\mathbf{240.0}$ \\
$(25, 0.9)$ & $0.007$ & $\mathbf{0.005}$ & $0.045$ & $0.049$ & $0.001$ & $\mathbf{0.001}$ & $227.0$ & $\mathbf{153.0}$ \\
$(25, 0.95)$ & $0.003$ & $\mathbf{0.003}$ & $0.039$ & $0.015$ & $0.000$ & $\mathbf{0.000}$ & $71.3$ & $\mathbf{46.3}$ \\
\cmidrule(lr){1-9}
$(30, 0.7)$ & $0.010$ & $\mathbf{0.008}$ & $6.409$ & $9.829$ & $\mathbf{0.041}$ & $0.053$ & $\mathbf{9825.3}$ & $10177.7$ \\
$(30, 0.75)$ & $0.010$ & $\mathbf{0.009}$ & $5.566$ & $5.855$ & $\mathbf{0.041}$ & $0.044$ & $9786.0$ & $\mathbf{8143.7}$ \\
$(30, 0.8)$ & $0.009$ & $\mathbf{0.006}$ & $1.780$ & $2.697$ & $\mathbf{0.026}$ & $0.030$ & $6220.0$ & $\mathbf{5869.3}$ \\
$(30, 0.85)$ & $0.017$ & $\mathbf{0.015}$ & $0.674$ & $0.805$ & $\mathbf{0.015}$ & $0.017$ & $3610.0$ & $\mathbf{3229.3}$ \\
$(30, 0.9)$ & $0.007$ & $\mathbf{0.005}$ & $0.161$ & $0.199$ & $\mathbf{0.005}$ & $0.007$ & $\mathbf{1232.0}$ & $1309.3$ \\
$(30, 0.95)$ & $0.007$ & $\mathbf{0.006}$ & $0.058$ & $0.076$ & $0.002$ & $\mathbf{0.002}$ & $380.7$ & $\mathbf{260.3}$ \\
\cmidrule(lr){1-9}
$(35, 0.7)$ & $0.024$ & $\mathbf{0.022}$ & $182.914$ & $142.446$ & $0.363$ & $\mathbf{0.311}$ & $72920.7$ & $\mathbf{49665.0}$ \\
$(35, 0.75)$ & $0.021$ & $\mathbf{0.019}$ & $112.531$ & $141.736$ & $\mathbf{0.277}$ & $0.306$ & $57875.7$ & $\mathbf{48878.7}$ \\
$(35, 0.8)$ & $0.013$ & $\mathbf{0.010}$ & $30.880$ & $36.765$ & $\mathbf{0.151}$ & $0.164$ & $31950.7$ & $\mathbf{27655.0}$ \\
$(35, 0.85)$ & $0.013$ & $\mathbf{0.010}$ & $6.115$ & $6.796$ & $\mathbf{0.061}$ & $0.069$ & $13374.0$ & $\mathbf{11582.0}$ \\
$(35, 0.9)$ & $0.013$ & $\mathbf{0.011}$ & $1.296$ & $0.799$ & $0.025$ & $\mathbf{0.021}$ & $5733.0$ & $\mathbf{3617.0}$ \\
$(35, 0.95)$ & $0.009$ & $\mathbf{0.007}$ & $0.359$ & $0.259$ & $0.006$ & $\mathbf{0.006}$ & $1199.3$ & $\mathbf{944.3}$ \\
\bottomrule
\end{tabular}

}
 \caption{Computational comparison on small bilinear instances. Results are averaged over three instances for each \((n,\gamma)\). Smallest values in each category are shown in bold.}

  \label{tb:bilinear_expt1}
\end{table}

The results are reported in Table~\ref{tb:bilinear_expt1}. Overall, \textsc{Lin} achieves the best runtime in almost all configurations, with \textsc{Vanilla} typically serving as the second-best method. This suggests that the compact linearized formulation is the most effective implementation for these instances. 
The two pure cut-based methods, \textsc{Cut} and \textsc{LiftCut}, are valid and solve nearly all instances to optimality, but their efficiency deteriorates as instances grow. In the most challenging configuration, \((n,\gamma)=(35,0.7)\), each method reaches the 300-second time limit in one of the random instances.
Their performance is also sensitive to the values of \(\gamma\). For example, when \(n=35\), \textsc{Cut} generates \(72920.7\) cuts at \(\gamma=0.7\), compared with only \(1199.3\) cuts at \(\gamma=0.95\). This is expected because a smaller threshold enlarges the feasible region and requires more cuts to exclude weak incumbent solutions.

Comparing the two cut-based variants, \textsc{Cut} generally has lower separation overhead, while \textsc{LiftCut} often produces stronger inequalities. For instance, at \((n,\gamma)=(35,0.7)\), \textsc{LiftCut} reduces the average number of cuts from \(72920.7\) to \(49665.0\); at \((n,\gamma)=(25,0.85)\), it reduces the count from \(428.0\) to \(240.0\). However, stronger cuts do not always yield faster solution times, since lifting adds separation effort and may interact differently with branch-and-bound. For example, at \((n,\gamma)=(35,0.75)\), \textsc{LiftCut} generates fewer cuts than \textsc{Cut}, but has a higher runtime. Thus, lifted cuts improve cut quality, but their computational benefit depends on whether the reduced cut count offsets the additional separation cost.

Overall, although the pure cut-based algorithms are valid, they are not as efficient as direct solution approaches for problem~\eqref{eq:bibilinear}. In Experiment~2, we therefore evaluate their impact when incorporated as root-node cuts for \textsc{Lin}.

\subsubsection{Root-Node Strengthening via Bimonotone and Lifted Cuts}

This experiment compares \textsc{Lin} with two root-cut variants that strengthen the root-node relaxation. 
\textsc{Lin+Cut} adds standard bimonotone cuts, while \textsc{Lin+LiftCut} adds lifted vitality cuts~\eqref{eq:vat}. 
At the root node, we sample $1000$ binary candidates by randomized rounding of the root LP solution and remove duplicates. 
For each infeasible candidate, we generate standard or lifted cuts using the corresponding separation routine, retain cuts violated by the root LP solution, remove duplicate cuts, and add up to $100$ cuts through Gurobi's \textsc{cbCut} in the root node callback. 
Cut generation is parallelized for efficiency.
We vary $n$ over $\{1000,1250,\ldots,2250\}$ and $\gamma$ over $\{0.7,0.75,\cdots,0.95\}$, with three random instances per configuration, yielding $108$ instances in total.
We impose a time limit of $600$ seconds for each run.

\vbadness=10001
\begin{table}[!tbp]
  \centering
\rev{
 \scriptsize
\setlength{\tabcolsep}{3pt}
\begin{tabular}{lrrr@{\hskip 10pt}rrr@{\hskip 10pt}rr}
\toprule
Config & \multicolumn{3}{c}{Solve Time (s)} & \multicolumn{3}{c}{Num. Nodes} & \multicolumn{2}{c}{Num. Root Cuts}\\
\cmidrule(lr){1-1}\cmidrule(lr){2-4}\cmidrule(lr){5-7}\cmidrule(lr){8-9}
$(n,\gamma)$ & Lin & Lin+Cut & Lin+LiftCut & Lin & Lin+Cut & Lin+LiftCut & Lin+Cut & Lin+LiftCut \\
\midrule
$(1000, 0.7)$ & $\mathbf{9.66}$ & $9.97$ & $10.20$ & $241$ & $241$ & $\mathbf{178}$ & $\mathbf{5.0}$ & $5.7$ \\
$(1000, 0.75)$ & $17.61$ & $17.94$ & $\mathbf{17.17}$ & $\mathbf{355}$ & $476$ & $370$ & $\mathbf{3.3}$ & $8.0$ \\
$(1000, 0.8)$ & $\mathbf{14.80}$ & $15.08$ & $15.19$ & $\mathbf{1259}$ & $\mathbf{1259}$ & $\mathbf{1259}$ & $\mathbf{3.3}$ & $4.7$ \\
$(1000, 0.85)$ & $\mathbf{23.36}$ & $23.40$ & $23.79$ & $845$ & $845$ & $\mathbf{841}$ & $\mathbf{3.7}$ & $7.7$ \\
$(1000, 0.9)$ & $14.35$ & $13.97$ & $\mathbf{13.93}$ & $28$ & $\mathbf{27}$ & $\mathbf{27}$ & $\mathbf{2.0}$ & $6.7$ \\
$(1000, 0.95)$ & $\mathbf{83.30}$ & $86.89$ & $84.96$ & $1390$ & $1440$ & $\mathbf{1340}$ & $\mathbf{0.7}$ & $1.7$ \\
\cmidrule(lr){1-9}
$(1250, 0.7)$ & $29.34$ & $\mathbf{29.23}$ & $29.26$ & $1109$ & $1217$ & $\mathbf{720}$ & $\mathbf{2.3}$ & $3.7$ \\
$(1250, 0.75)$ & $\mathbf{22.56}$ & $25.92$ & $27.38$ & $\mathbf{1199}$ & $1386$ & $1358$ & $\mathbf{3.3}$ & $5.3$ \\
$(1250, 0.8)$ & $\mathbf{30.62}$ & $31.77$ & $31.41$ & $\mathbf{965}$ & $\mathbf{965}$ & $1496$ & $\mathbf{2.3}$ & $3.3$ \\
$(1250, 0.85)$ & $\mathbf{69.56}$ & $72.72$ & $74.71$ & $4437$ & $\mathbf{4150}$ & $4380$ & $\mathbf{5.0}$ & $8.7$ \\
$(1250, 0.9)$ & $88.34$ & $\mathbf{76.56}$ & $77.21$ & $1695$ & $\mathbf{1665}$ & $\mathbf{1665}$ & $\mathbf{1.0}$ & $4.7$ \\
$(1250, 0.95)$ & $\mathbf{153.71}$ & $158.40$ & $158.39$ & $\mathbf{945}$ & $\mathbf{945}$ & $\mathbf{945}$ & $\mathbf{0.0}$ & $\mathbf{0.0}$ \\
\cmidrule(lr){1-9}
$(1500, 0.7)$ & $32.48$ & $\mathbf{32.42}$ & $33.97$ & $1442$ & $794$ & $\mathbf{677}$ & $\mathbf{6.7}$ & $9.3$ \\
$(1500, 0.75)$ & $29.32$ & $30.87$ & $\mathbf{28.47}$ & $\mathbf{1553}$ & $1586$ & $1728$ & $\mathbf{4.7}$ & $7.0$ \\
$(1500, 0.8)$ & $\mathbf{116.24}$ & $120.06$ & $149.73$ & $\mathbf{3908}$ & $\mathbf{3908}$ & $4418$ & $\mathbf{6.0}$ & $11.3$ \\
$(1500, 0.85)$ & $113.70$ & $114.40$ & $\mathbf{63.38}$ & $4084$ & $4084$ & $\mathbf{2192}$ & $\mathbf{2.7}$ & $6.0$ \\
$(1500, 0.9)$ & $35.39$ & $\mathbf{35.37}$ & $35.70$ & $1541$ & $1556$ & $\mathbf{1139}$ & $\mathbf{1.3}$ & $4.3$ \\
$(1500, 0.95)$ & $600.46$ & $\mathbf{503.70}$ & $600.49$ & $1299$ & $\mathbf{1151}$ & $1255$ & $\mathbf{0.3}$ & $1.7$ \\
\cmidrule(lr){1-9}
$(1750, 0.7)$ & $\mathbf{51.07}$ & $51.74$ & $55.05$ & $6612$ & $6612$ & $\mathbf{5749}$ & $\mathbf{5.3}$ & $6.7$ \\
$(1750, 0.75)$ & $\mathbf{64.86}$ & $76.96$ & $94.23$ & $2068$ & $\mathbf{1273}$ & $2498$ & $\mathbf{3.7}$ & $7.0$ \\
$(1750, 0.8)$ & $98.58$ & $\mathbf{88.30}$ & $88.57$ & $\mathbf{1569}$ & $2868$ & $4690$ & $\mathbf{5.0}$ & $11.0$ \\
$(1750, 0.85)$ & $\mathbf{83.45}$ & $85.49$ & $89.41$ & $\mathbf{320}$ & $321$ & $476$ & $\mathbf{4.0}$ & $6.3$ \\
$(1750, 0.9)$ & $232.71$ & $\mathbf{64.93}$ & $65.38$ & $2729$ & $3110$ & $\mathbf{2497}$ & $\mathbf{1.3}$ & $4.3$ \\
$(1750, 0.95)$ & $\mathbf{515.00}$ & $517.34$ & $518.86$ & $977$ & $1038$ & $\mathbf{918}$ & $\mathbf{1.0}$ & $1.7$ \\
\cmidrule(lr){1-9}
$(2000, 0.7)$ & $\mathbf{64.41}$ & $65.74$ & $74.89$ & $\mathbf{1271}$ & $1679$ & $1379$ & $\mathbf{2.7}$ & $4.3$ \\
$(2000, 0.75)$ & $89.55$ & $\mathbf{89.36}$ & $137.79$ & $\mathbf{1841}$ & $\mathbf{1841}$ & $3408$ & $\mathbf{2.0}$ & $3.3$ \\
$(2000, 0.8)$ & $\mathbf{90.58}$ & $90.70$ & $90.82$ & $2315$ & $2344$ & $\mathbf{1095}$ & $\mathbf{2.0}$ & $5.3$ \\
$(2000, 0.85)$ & $168.95$ & $\mathbf{161.87}$ & $166.90$ & $2407$ & $2407$ & $\mathbf{1802}$ & $\mathbf{5.3}$ & $7.0$ \\
$(2000, 0.9)$ & $423.58$ & $255.96$ & $\mathbf{254.99}$ & $\mathbf{721}$ & $1701$ & $1702$ & $\mathbf{1.7}$ & $5.3$ \\
$(2000, 0.95)$ & $\mathbf{600.55}$ & $600.66$ & $600.67$ & $492$ & $479$ & $\mathbf{462}$ & $\mathbf{0.0}$ & $\mathbf{0.0}$ \\
\cmidrule(lr){1-9}
$(2250, 0.7)$ & $69.05$ & $68.36$ & $\mathbf{63.92}$ & $\mathbf{2169}$ & $\mathbf{2169}$ & $2196$ & $\mathbf{4.0}$ & $8.0$ \\
$(2250, 0.75)$ & $130.17$ & $131.02$ & $\mathbf{129.59}$ & $2542$ & $2542$ & $\mathbf{1535}$ & $\mathbf{4.0}$ & $6.7$ \\
$(2250, 0.8)$ & $180.41$ & $184.20$ & $\mathbf{165.65}$ & $\mathbf{2078}$ & $2609$ & $2798$ & $\mathbf{3.0}$ & $4.7$ \\
$(2250, 0.85)$ & $\mathbf{98.00}$ & $102.11$ & $98.80$ & $\mathbf{2818}$ & $3133$ & $3530$ & $\mathbf{1.7}$ & $4.3$ \\
$(2250, 0.9)$ & $\mathbf{433.21}$ & $433.40$ & $433.66$ & $1580$ & $\mathbf{1573}$ & $1676$ & $\mathbf{0.7}$ & $3.7$ \\
$(2250, 0.95)$ & $600.83$ & $600.91$ & $\mathbf{600.59}$ & $\mathbf{615}$ & $807$ & $622$ & $\mathbf{1.7}$ & $2.3$ \\
\bottomrule
\end{tabular}

}
 \caption{Computational comparison on large bilinear instances. Results are averaged over three instances for each \((n,\gamma)\). Smallest values in each category are shown in bold.}
  \label{tb:bilinear_expt2}
\end{table}

The results are reported in Table~\ref{tb:bilinear_expt2}. The impact of adding root-node cuts to \textsc{Lin} is mixed but often beneficial. In many configurations, either \textsc{Lin+Cut} or \textsc{Lin+LiftCut} improves the solve time, especially on larger or harder instances. For example, at \((n,\gamma)=(1750,0.9)\), the runtime decreases from \(232.71\) seconds for \textsc{Lin} to about \(65\) seconds with root-node cuts. Similarly, at \((n,\gamma)=(2000,0.9)\), the runtime decreases from \(423.58\) seconds to about \(255\) seconds. However, the benefit is not uniform: on easier instances, the extra cost of generating and processing root cuts can offset the tightening effect.

Root-node cuts often reduce the branch-and-bound tree, suggesting that they strengthen the root relaxation. For $20$ out of the \(36\) configurations, the cut-based variants lead to fewer explored nodes. Moreover, \textsc{Lin+LiftCut} often produces stronger root-level tightening than \textsc{Lin+Cut}. This is consistent with the larger number of root cuts generated by \textsc{Lin+LiftCut}, since lifted cuts can be derived from multiple enlarged sets in Example~\ref{eg:bilift}.

Overall, these results show that the bimonotone cuts and their lifted variants are not uniformly superior as root-node enhancements, but they can provide meaningful performance gains on selected larger and harder bilinear instances. 
The proposed framework helps identify and generate diverse families of cuts that can be used for this root-tightening purpose.

\end{revision}



\subsection{\rev{Case Study 2: Distributionally Robust Network Site Selection}}
\begin{revision}
  Chance-constrained binary optimization provides a natural framework for selecting discrete resources under uncertain feasibility requirements. This class of models has been studied through chance-constrained set packing~\citep{song2014chance} and probabilistic set covering problems~\citep{beraldi2002probabilistic,saxena2010mip}, with applications in vaccine allocation~\citep{tanner2010iis} and capital rationing~\citep{beraldi2012capital}.
Our case study extends this to a distributionally robust optimization (DRO) that hedges against distributional shift~\citep{mohajerin2018data,gao2023distributionally}, demonstrating how interval system identification enables multiple cut-based solution strategies within the proposed framework.
\end{revision}


Consider a supply network $G = (V, E)$ where a company aims to select certain nodes from $V$ to supply products to all demand nodes. If constructed, each site $i \in V$ would induce a fixed net benefit $r_i$, which incorporates both economic and environmental factors, and thus can be either positive or negative. The potential supplies from all nodes to node $j$ are denoted by the vector $a_j = (a_{ij})_{i \in V} \in \mathbb{R}_+^{|V|}$, and the demand at node $j$ is $b_j \in \mathbb{R}_+$. Historical data indicates that, for every $j \in V$, $(a_j, b_j)$ is a random vector following some \rev{discrete} empirical joint distribution $\bar{\mathbb{P}}_j$ with the \rev{finite-scenario} support denoted by $\Xi_j$. 
\rev{
To hedge against potential distributional inaccuracy or shift, we define the Wasserstein ambiguity set
\[
\mathfrak P(\bar{\mathbb P}_j)
:=
\{\mathbb P \mid W_1(\mathbb P,\bar{\mathbb P}_j)\leq \eta\}
\text{ with }
W_1(\mathbb P,\bar{\mathbb P}_j)
:=
\inf_{\pi\in\Pi(\mathbb P,\bar{\mathbb P}_j)}
\mathbb E_{\pi}\!\left[
\left\|(a_j,b_j)-(a'_j,b'_j)\right\|
\right],
\]
where \(\eta\) is the ambiguity radius, \(W_1\) denotes the Wasserstein-1 distance, and \(\|\cdot\|\) is chosen as the 2-norm.
}
To avoid service overlap, the company also requires that the selected sites $T \subseteq V$ form an independent set. Then, the following formulation uses chance constraints to provide demand satisfaction guarantee.
  \begin{subequations}
    \label{eq:case}
  \begin{align} 
    \max\limits_{x \in \{0,1\}^V} &~ \sum_{i \in V} r_i x_i \label{eq:case00}\\
    \text{s.t.} & ~\sup_{\mathbb P \in \mathfrak P(\bar{\mathbb P}_j)}\mathbb P\left(\sum_{i \in \delta[j]}a_{ij}x_i < b_j\right)\leq \epsilon_j, ~\forall j \in V\label{eq:case01}\\
                & ~ \sum_{i \in C} x_i \leq 1, ~~~\quad\quad\quad\quad\quad\quad\quad\quad\forall C=\{i,j\} \in E.\label{eq:case02}
\end{align}
\end{subequations}
The objective \eqref{eq:case00} is to maximize the net benefit. \eqref{eq:case01} is the \rev{DRO} chance constraint to impose a demand satisfaction level of $1-\epsilon_j$ with $\delta[j]$ denoting the closed neighborhood of vertex $j$. Constraint \eqref{eq:case02} ensures that the resulting solution is an independent set, which can be further enhanced using clique elimination constraints.

From the perspective of Section~\ref{sec:cgalgebra}, constraints \eqref{eq:case01} and \eqref{eq:case02} define two monotone subsystems: $\Omega_{\mathcal{X}_1}$ is an upper system (\rev{DRO} chance constraints with nonnegative coefficients), while $\Omega_{\mathcal{X}_2}$ is a lower system (independent set constraints). Hence, the feasible region $\Omega=\Omega_{\mathcal{X}_1}\cap\Omega_{\mathcal{X}_2}$ forms an interval system. This classification yields multiple reformulation paths: covering inequalities can be applied to $\Omega_{\mathcal{X}_1}$, elimination inequalities to $\Omega_{\mathcal{X}_2}$, or both simultaneously. In particular, the \rev{DRO chance-constraint \eqref{eq:case01} can be exactly reformulated into the following inequalities by Theorem~\ref{thm:reform}}
\begin{equation}
\label{eq:sites}
\sum_{i\in T} x_i \ge 1,\qquad \forall T \in m\ec{\nc{(\Omega_{\mathcal X_1})}}.
\end{equation}
A direct computation shows that each $T$ in this index set corresponds to a minimal collection of sites whose complement yields a satisfaction probability strictly less than $1-\epsilon_j$. Theorem~\ref{thm:reform} guarantees that this reformulation is exact, while Proposition~\ref{prop:facet} characterizes the facet-defining inequalities. 
This illustrates how the monotone system perspective not only enables flexible reformulations but also unifies their analysis and strengthening within a broader polyhedral toolkit.

\subsubsection{Experiment Setup}
Since the proposed method supports cut generation for arbitrary binary spaces or subspaces, it enables flexible algorithmic choices. We compare the following four implementations, each of which corresponds to a different way of exploiting the two monotone systems $\Omega_{\mathcal X_1}$ and $\Omega_{\mathcal X_2}$.
\begin{itemize}
  \item \textsc{NoCut}: Baseline finite scenario approximation (FSA) algorithm; no monotone cuts used.
  \item \textsc{ClqCut}: Exploits lower system $\Omega_{\mathcal X_2}$; adds clique-based elimination cuts.
  \item \textsc{SatCut}: Exploits upper system $\Omega_{\mathcal X_1}$; generates covering cuts \eqref{eq:sites}.
  \item \textsc{AllCut}: Uses both, consistent with interval-system structure $\Omega_{\mathcal X_1} \cap \Omega_{\mathcal X_2}$.
\end{itemize}
\begin{revision}
  In \textsc{NoCut} and \textsc{ClqCut}, we rewrite \eqref{eq:case01} as follows
\[
  \sup_{\mathbb P_j \in \mathfrak P(\bar{\mathbb P}_j)} \mathbb P_j\left(\sum_{i \in \delta[j]} a_{ij}x_i < b_j\right) = \sup_{\mathbb P_j \in \mathfrak P(\bar{\mathbb P}_j)} \mathbb E_{\mathbb P_j}\left[\mathbb I_{\Xi_j(x)}\right]\leq \epsilon_j , ~\forall j \in V,
\]
where $\Xi_j(x):=\{(a_{j}, b_j) \in \Xi_j \mid \sum_{i \in \delta[j]}a_{ij}x_i < b_j\}$, $\mathbb I_{\Xi_j(x)}$ is its indicator function.
Applying standard finite-support DRO reformulation~\citep{mohajerin2018data,luo2020distributionally}, we obtain  
\begin{equation}
  \label{eq:dro}
\sup_{\mathbb{P}_j \in \mathfrak{P}(\bar{\mathbb{P}}_j)} \mathbb{E}_{\mathbb{P}_j}\!\left[\mathbb{I}_{\Xi_j(x)}\right]
\;=\;
\left\{
\begin{aligned}
\inf_{\gamma_j\geq 0,\, s_j^k} \quad & \sum_{k \in [K]} p^k_j s_j^k + \eta\, \gamma_j \\
\text{s.t.} \quad & s_j^k + \gamma_j \big\| (a_j^l, b_j^l) - (a_j^k, b_j^k) \big\| \geq \mathbb{I}_{\Xi_j(x)}(a_j^l, b_j^l), ~ \forall\, k, l \in [K],
\end{aligned}
\right.
\end{equation}
where \(p_j^k\) is the nominal probability of scenario \(k\) for node \(j\). 
Let $z^k_j \in \{0,1\}$ indicate constraint satisfaction under the scenario $k$. We obtain the following constraint representation of~\eqref{eq:case01}.
$$
\begin{aligned}
& \sum_{k \in [K]} p_j^k s_j^k + \eta\, \gamma_j  \leq \epsilon_j, ~\forall j \in V\\
& s_j^k + \gamma_j \big\| (a_j^l, b_j^l) - (a_j^k, b_j^k) \big\| \geq 1- z^l_j, ~ \forall\, k, l \in [K]\\
&  \sum_{i \in \delta[j]} a^k_{ij} x_i  + b_j^k(1-z^k_j) \geq b^k_j,~\forall k \in [K], j \in V.
\end{aligned}
$$
\end{revision}
In \textsc{ClqCut}, we generate cuts for \eqref{eq:case02} on-the-fly by identifying up to three violated \emph{maximal cliques} (with $|C|\ge 2$) in the induced subgraph $G[T]$ at an incumbent solution $x_T$ from the master problem. In \textsc{SatCut}, given such an $x_T$, the $j$th subproblem \rev{solves~\eqref{eq:dro} to evaluate the worst-case violation} probability. This will be used to either confirm or reject the feasibility of $x_T$ regarding \eqref{eq:case01}. For the rejection case, the proposed constraint \eqref{eq:sites} will be added to the master according to \eqref{eq:decomp}. \rev{In \textsc{AllCut}, both types of cuts are separated on-the-fly.}

We generate \(a_{ij}\) from \(\mathrm{Beta}(2,2)\) if \(\{i,j\}\in E\), and set \(a_{ij}=0\) otherwise. The demands \(b_j\) are drawn uniformly from \([0,0.1]\), and the benefit coefficients \(r_i\) are integers drawn uniformly from \(\{-20,\dots,20\}\).
The experimental parameters are \(n\in\{40,80,120\}\), \(m\in\{0.3n(n-1)/2,0.7n(n-1)/2\}\), \(\epsilon\in\{0.05,0.1\}\), and \(K\in\{200,500,800\}\), representing the number of vertices, number of edges, demand violation tolerance, and number of scenarios, respectively. Each tuple \((n,m,\epsilon,K)\) defines an instance configuration, where the two choices of \(m\) correspond to graph densities \(0.3\) and \(0.7\). For each configuration, we generate three connected Erd\H{o}s--R\'{e}nyi graphs, resulting in \(108\) instances in total. All four algorithms are run on each instance with a \(600\)-second time limit. 
\subsubsection{Results Analysis}

\vbadness=10001
\begin{table}[!tbp]
  \centering
\rev{
 \scriptsize
\begin{tabular}{lrrrrrrrrrr}
\toprule
Config & \multicolumn{4}{c}{Runtime (s)} & \multicolumn{3}{c}{Cuts Sep. Time (s)} & \multicolumn{3}{c}{Num. of Cuts}\\
\cmidrule(lr){1-1}\cmidrule(lr){2-5}\cmidrule(lr){6-8}\cmidrule(lr){9-11}
\((n,m,\epsilon,K)\) & NoCut & ClqCut & SatCut & AllCut & ClqCut & SatCut & AllCut & ClqCut & SatCut & AllCut \\
\midrule
\((40, 234, 0.05, 200)\) & \(50.49\) & \(91.37\) & \(\mathbf{0.55}\) & \(0.79\) & \(\mathbf{0.01}\) & \(0.47\) & \(0.78\) & \(64.7\) & \(\mathbf{36.7}\) & \(87.0\) \\
\((40, 234, 0.05, 500)\) & \(146.26\) & \(474.70\) & \(\mathbf{0.87}\) & \(1.45\) & \(\mathbf{0.01}\) & \(0.81\) & \(1.45\) & \(52.0\) & \(\mathbf{12.0}\) & \(57.0\) \\
\((40, 234, 0.05, 800)\) & \(365.12\) & \(603.64\) & \(\mathbf{5.54}\) & \(10.01\) & \(\mathbf{0.00}\) & \(5.40\) & \(9.99\) & \(36.7\) & \(\mathbf{21.3}\) & \(98.7\) \\
\cmidrule(lr){2-11}
\((40, 234, 0.1, 200)\) & \(17.79\) & \(52.60\) & \(\mathbf{0.25}\) & \(0.40\) & \(\mathbf{0.00}\) & \(0.21\) & \(0.40\) & \(44.0\) & \(\mathbf{21.7}\) & \(50.7\) \\
\((40, 234, 0.1, 500)\) & \(163.50\) & \(459.72\) & \(\mathbf{1.63}\) & \(1.83\) & \(\mathbf{0.00}\) & \(1.52\) & \(1.83\) & \(48.3\) & \(\mathbf{31.3}\) & \(76.7\) \\
\((40, 234, 0.1, 800)\) & \(530.77\) & \(604.36\) & \(\mathbf{2.81}\) & \(5.84\) & \(\mathbf{0.00}\) & \(2.71\) & \(5.83\) & \(32.7\) & \(\mathbf{11.7}\) & \(64.0\) \\
\midrule
\((40, 546, 0.05, 200)\) & \(15.57\) & \(50.28\) & \(\mathbf{0.18}\) & \(0.52\) & \(\mathbf{0.01}\) & \(0.17\) & \(0.52\) & \(56.3\) & \(\mathbf{5.0}\) & \(48.0\) \\
\((40, 546, 0.05, 500)\) & \(127.26\) & \(222.93\) & \(\mathbf{0.22}\) & \(0.98\) & \(\mathbf{0.01}\) & \(0.21\) & \(0.97\) & \(49.3\) & \(\mathbf{0.0}\) & \(44.7\) \\
\((40, 546, 0.05, 800)\) & \(225.96\) & \(450.42\) & \(\mathbf{2.63}\) & \(3.49\) & \(\mathbf{0.01}\) & \(2.59\) & \(3.48\) & \(46.0\) & \(\mathbf{7.3}\) & \(57.3\) \\
\cmidrule(lr){2-11}
\((40, 546, 0.1, 200)\) & \(19.19\) & \(39.39\) & \(\mathbf{0.09}\) & \(0.46\) & \(\mathbf{0.00}\) & \(0.08\) & \(0.45\) & \(44.0\) & \(\mathbf{3.7}\) & \(52.7\) \\
\((40, 546, 0.1, 500)\) & \(138.13\) & \(365.83\) & \(\mathbf{0.44}\) & \(1.20\) & \(\mathbf{0.01}\) & \(0.42\) & \(1.19\) & \(51.0\) & \(\mathbf{3.7}\) & \(49.3\) \\
\((40, 546, 0.1, 800)\) & \(386.62\) & \(501.88\) & \(\mathbf{2.76}\) & \(3.23\) & \(\mathbf{0.01}\) & \(2.69\) & \(3.23\) & \(45.7\) & \(\mathbf{10.3}\) & \(55.3\) \\
\midrule
\((80, 948, 0.05, 200)\) & \(118.98\) & \(208.08\) & \(\mathbf{0.60}\) & \(2.25\) & \(\mathbf{0.02}\) & \(0.45\) & \(2.22\) & \(160.7\) & \(\mathbf{20.0}\) & \(179.7\) \\
\((80, 948, 0.05, 500)\) & \(381.77\) & \(603.82\) & \(\mathbf{0.92}\) & \(5.34\) & \(\mathbf{0.02}\) & \(0.82\) & \(5.32\) & \(74.0\) & \(\mathbf{12.0}\) & \(114.3\) \\
\((80, 948, 0.05, 800)\) & \(613.39\) & \(611.07\) & \(\mathbf{8.99}\) & \(18.90\) & \(\mathbf{0.04}\) & \(8.67\) & \(18.87\) & \(46.0\) & \(\mathbf{27.7}\) & \(156.0\) \\
\cmidrule(lr){2-11}
\((80, 948, 0.1, 200)\) & \(126.13\) & \(215.67\) & \(\mathbf{1.80}\) & \(3.01\) & \(\mathbf{0.02}\) & \(1.35\) & \(2.98\) & \(184.0\) & \(\mathbf{87.0}\) & \(257.0\) \\
\((80, 948, 0.1, 500)\) & \(532.89\) & \(606.79\) & \(\mathbf{2.78}\) & \(9.19\) & \(\mathbf{0.02}\) & \(2.52\) & \(9.16\) & \(77.0\) & \(\mathbf{23.7}\) & \(179.7\) \\
\((80, 948, 0.1, 800)\) & \(614.56\) & \(607.52\) & \(\mathbf{6.78}\) & \(22.53\) & \(\mathbf{0.03}\) & \(6.34\) & \(22.50\) & \(54.0\) & \(\mathbf{16.7}\) & \(169.3\) \\
\midrule
\((80, 2212, 0.05, 200)\) & \(144.11\) & \(183.95\) & \(\mathbf{0.41}\) & \(2.39\) & \(\mathbf{0.02}\) & \(0.30\) & \(2.36\) & \(197.3\) & \(\mathbf{9.3}\) & \(155.0\) \\
\((80, 2212, 0.05, 500)\) & \(602.93\) & \(609.06\) & \(\mathbf{1.79}\) & \(8.07\) & \(\mathbf{0.02}\) & \(1.61\) & \(8.04\) & \(69.3\) & \(\mathbf{14.7}\) & \(190.0\) \\
\((80, 2212, 0.05, 800)\) & \(612.94\) & \(636.94\) & \(\mathbf{2.87}\) & \(10.76\) & \(\mathbf{0.03}\) & \(2.68\) & \(10.73\) & \(43.7\) & \(\mathbf{5.7}\) & \(126.7\) \\
\cmidrule(lr){2-11}
\((80, 2212, 0.1, 200)\) & \(118.73\) & \(158.51\) & \(\mathbf{0.48}\) & \(2.60\) & \(\mathbf{0.02}\) & \(0.36\) & \(2.57\) & \(205.0\) & \(\mathbf{13.0}\) & \(167.0\) \\
\((80, 2212, 0.1, 500)\) & \(602.96\) & \(610.77\) & \(\mathbf{2.36}\) & \(7.46\) & \(\mathbf{0.02}\) & \(2.17\) & \(7.43\) & \(61.0\) & \(\mathbf{13.3}\) & \(174.7\) \\
\((80, 2212, 0.1, 800)\) & \(611.52\) & \(657.86\) & \(\mathbf{3.19}\) & \(12.42\) & \(\mathbf{0.05}\) & \(2.94\) & \(12.38\) & \(45.7\) & \(\mathbf{8.0}\) & \(159.0\) \\
\midrule
\((120, 2142, 0.05, 200)\) & \(249.54\) & \(584.33\) & \(\mathbf{1.60}\) & \(7.82\) & \(\mathbf{0.09}\) & \(0.85\) & \(7.74\) & \(423.0\) & \(\mathbf{27.7}\) & \(443.3\) \\
\((120, 2142, 0.05, 500)\) & \(605.35\) & \(607.03\) & \(\mathbf{5.95}\) & \(30.82\) & \(\mathbf{0.03}\) & \(4.95\) & \(30.71\) & \(83.0\) & \(\mathbf{47.3}\) & \(472.7\) \\
\((120, 2142, 0.05, 800)\) & \(629.46\) & \(617.93\) & \(\mathbf{5.68}\) & \(32.03\) & \(\mathbf{0.07}\) & \(5.00\) & \(31.98\) & \(52.0\) & \(\mathbf{17.3}\) & \(280.3\) \\
\cmidrule(lr){2-11}
\((120, 2142, 0.1, 200)\) & \(279.09\) & \(493.74\) & \(\mathbf{1.88}\) & \(6.81\) & \(\mathbf{0.06}\) & \(1.06\) & \(6.74\) & \(404.3\) & \(\mathbf{40.0}\) & \(414.0\) \\
\((120, 2142, 0.1, 500)\) & \(605.85\) & \(604.74\) & \(\mathbf{4.96}\) & \(23.13\) & \(\mathbf{0.03}\) & \(3.94\) & \(23.05\) & \(86.0\) & \(\mathbf{25.7}\) & \(395.3\) \\
\((120, 2142, 0.1, 800)\) & \(654.43\) & \(615.76\) & \(\mathbf{13.64}\) & \(41.70\) & \(\mathbf{0.09}\) & \(12.49\) & \(41.63\) & \(57.0\) & \(\mathbf{34.3}\) & \(353.3\) \\
\midrule
\((120, 4998, 0.05, 200)\) & \(342.62\) & \(521.95\) & \(\mathbf{0.96}\) & \(5.90\) & \(\mathbf{0.07}\) & \(0.41\) & \(5.85\) & \(482.0\) & \(\mathbf{5.7}\) & \(297.0\) \\
\((120, 4998, 0.05, 500)\) & \(606.07\) & \(603.81\) & \(\mathbf{2.25}\) & \(20.04\) & \(\mathbf{0.03}\) & \(1.50\) & \(19.99\) & \(61.3\) & \(\mathbf{9.0}\) & \(368.7\) \\
\((120, 4998, 0.05, 800)\) & \(619.42\) & \(614.30\) & \(\mathbf{3.46}\) & \(21.09\) & \(\mathbf{0.08}\) & \(2.98\) & \(21.05\) & \(28.0\) & \(\mathbf{2.3}\) & \(233.3\) \\
\cmidrule(lr){2-11}
\((120, 4998, 0.1, 200)\) & \(181.65\) & \(416.98\) & \(\mathbf{0.61}\) & \(4.33\) & \(\mathbf{0.05}\) & \(0.42\) & \(4.31\) & \(331.0\) & \(\mathbf{7.3}\) & \(223.7\) \\
\((120, 4998, 0.1, 500)\) & \(606.97\) & \(605.07\) & \(\mathbf{1.49}\) & \(15.11\) & \(\mathbf{0.03}\) & \(0.79\) & \(15.06\) & \(60.0\) & \(\mathbf{5.0}\) & \(295.3\) \\
\((120, 4998, 0.1, 800)\) & \(614.64\) & \(615.93\) & \(\mathbf{3.97}\) & \(25.06\) & \(\mathbf{0.07}\) & \(2.99\) & \(25.02\) & \(29.0\) & \(\mathbf{8.3}\) & \(260.0\) \\
\bottomrule
\end{tabular}

}
 \caption{Computational comparison for the site selection problem with DRO chance constraints. Results are averaged over three instances for each \((n,m,\epsilon,K)\). Smallest values in each category are shown in bold.}
  \label{tb:expt1}
\end{table}

\begin{revision}
The results are summarized in Table~\ref{tb:expt1}. Overall, \textsc{SatCut} substantially outperforms the other implementations across all configurations, often reducing solution times from hundreds of seconds to only a few seconds. This indicates that the satisfaction cuts derived from the upper-system structure of the DRO chance constraints provide the most effective tightening in this experiment. In contrast, \textsc{NoCut} and \textsc{ClqCut} frequently approach or reach the time limit as \(n\), graph density, or the number of scenarios \(K\) increases. The clique cuts in \textsc{ClqCut} do not consistently improve performance over \textsc{NoCut}, suggesting that tightening only the independent set constraints is less effective for these instances.

The comparison between \textsc{SatCut} and \textsc{AllCut} further shows that adding more cuts is not always beneficial. Although \textsc{AllCut} also exploits the satisfaction-cut structure, it typically generates many more cuts and spends more time in separation, which makes it slower than \textsc{SatCut}. The ``Num. of Cuts'' column shows that \textsc{SatCut} achieves its speedup with relatively few cuts, indicating that these cuts are strong and targeted. As expected, larger instances, denser graphs, and larger \(K\) generally increase computational difficulty, especially for \textsc{NoCut} and \textsc{ClqCut}. In contrast, the violation tolerance \(\epsilon\) has a less systematic effect on runtime. 
These results show that the proposed framework enables flexible and effective cut generation strategies by identifying useful monotone and interval subsystems. 
\end{revision}

\begin{revision}
\subsection{Summary of Numerical Insights}
Although Theorem~\ref{thm:reform} establishes that monotone solution spaces can be represented exactly using covering and elimination inequalities, the computational results suggest a broader utility for the proposed framework. When the underlying structure already admits a compact alternative formulation, such as the McCormick linearization of the bilinear terms in Case Study~1 or the independent set constraints in Case Study~2, full cut-based separation may increase solution time. In such cases, selectively adding strengthened inequalities, especially as root-node cuts, can still be beneficial on larger or harder instances. In contrast, when the monotone system involves a complex constraint representation, such as the DRO chance constraints in Case Study~2, a cut-based master--subproblem implementation can substantially reduce solution time. Overall, the proposed framework helps identify monotone, bimonotone, and interval structures from which valid and lifted inequalities can be systematically generated, while the most effective use of these inequalities should be tailored to the structure of each problem.
\end{revision}


\section{Conclusion}
\label{sec:conclusion}

\begin{revision}
This paper develops a unified framework for set-system approximation in binary integer programs (BIPs), interpreting covering and elimination inequalities as primitives for characterizing and approximating binary solution spaces. We show that these inequalities exactly describe the tightest monotone inner and outer approximations of arbitrary set systems, and become exact reformulations precisely under monotonicity. Building on this perspective, we recover classical structural correspondences, extend polyhedral tools from set covering to broader BIP settings, and develop reformulation methods for nonlinear and latent monotone structures, including auxiliary-variable-free bilinear reformulations, bimonotone cuts, and interval system decompositions.

The computational experiments illustrate that these reformulations have diverse uses beyond standalone exact models. Depending on the structure of the underlying monotone subsystem, the resulting inequalities can serve as root-node strengthening cuts, support master--subproblem cut generation, or guide hybrid reformulations. Future work may explore additional set system operators, scalable separation strategies, and extensions to mixed integer nonlinear programs.
\end{revision}

\newpage
\bibliographystyle{abbrvnat}
\bibliography{bibi/myref}

\begin{thebibliography}{71}
\providecommand{\natexlab}[1]{#1}
\providecommand{\url}[1]{\texttt{#1}}
\expandafter\ifx\csname urlstyle\endcsname\relax
  \providecommand{\doi}[1]{doi: #1}\else
  \providecommand{\doi}{doi: \begingroup \urlstyle{rm}\Url}\fi

\bibitem[Achuthan et~al.(1996)Achuthan, Caccetta, and Hill]{achuthan1996new}
N.~Achuthan, L.~Caccetta, and S.~Hill.
\newblock A new subtour elimination constraint for the vehicle routing problem.
\newblock \emph{European Journal of Operational Research}, 91\penalty0 (3):\penalty0 573--586, 1996.

\bibitem[Balas(1979)]{balas1979disjunctive}
E.~Balas.
\newblock Disjunctive programming.
\newblock \emph{Annals of discrete mathematics}, 5:\penalty0 3--51, 1979.

\bibitem[Balas and Ng(1989{\natexlab{a}})]{balas1989set}
E.~Balas and S.~M. Ng.
\newblock On the set covering polytope: I. all the facets with coefficients in $\{$0, 1, 2$\}$.
\newblock \emph{Mathematical Programming}, 43\penalty0 (1):\penalty0 57--69, 1989{\natexlab{a}}.

\bibitem[Balas and Ng(1989{\natexlab{b}})]{balas1989setii}
E.~Balas and S.~M. Ng.
\newblock On the set covering polytope: Ii. lifting the facets with coefficients in $\{$0, 1, 2$\}$.
\newblock \emph{Mathematical Programming}, 45\penalty0 (1):\penalty0 1--20, 1989{\natexlab{b}}.

\bibitem[Balas and Padberg(1972)]{balas1972set}
E.~Balas and M.~W. Padberg.
\newblock On the set-covering problem.
\newblock \emph{Operations Research}, 20\penalty0 (6):\penalty0 1152--1161, 1972.

\bibitem[Barahona and Mahjoub(1986)]{barahona1986cut}
F.~Barahona and A.~R. Mahjoub.
\newblock On the cut polytope.
\newblock \emph{Mathematical programming}, 36\penalty0 (2):\penalty0 157--173, 1986.

\bibitem[Barahona and Mahjoub(1989)]{barahona1989facets}
F.~Barahona and A.~R. Mahjoub.
\newblock Facets of the balanced (acyclic) induced subgraph polytope.
\newblock \emph{Mathematical programming}, 45\penalty0 (1):\penalty0 21--33, 1989.

\bibitem[Barahona et~al.(1985)Barahona, Gr{\"o}tschel, and Mahjoub]{barahona1985facets}
F.~Barahona, M.~Gr{\"o}tschel, and A.~R. Mahjoub.
\newblock Facets of the bipartite subgraph polytope.
\newblock \emph{Mathematics of Operations Research}, 10\penalty0 (2):\penalty0 340--358, 1985.

\bibitem[Baycik and Sullivan(2019)]{baycik2019robust}
N.~O. Baycik and K.~M. Sullivan.
\newblock Robust location of hidden interdictions on a shortest path network.
\newblock \emph{IISE Transactions}, 51\penalty0 (12):\penalty0 1332--1347, 2019.

\bibitem[Bellmore and Nemhauser(1968)]{bellmore1968traveling}
M.~Bellmore and G.~L. Nemhauser.
\newblock The traveling salesman problem: a survey.
\newblock \emph{Operations Research}, 16\penalty0 (3):\penalty0 538--558, 1968.

\bibitem[Benders(1962)]{bnnobrs1962partitioning}
J.~Benders.
\newblock Partitioning procedures for solving mixed-variables programming problems.
\newblock \emph{Numer. Math}, 4\penalty0 (1):\penalty0 238--252, 1962.

\bibitem[Beraldi and Ruszczy{\'n}ski(2002)]{beraldi2002probabilistic}
P.~Beraldi and A.~Ruszczy{\'n}ski.
\newblock The probabilistic set-covering problem.
\newblock \emph{Operations Research}, 50\penalty0 (6):\penalty0 956--967, 2002.

\bibitem[Beraldi et~al.(2012)Beraldi, Bruni, and Violi]{beraldi2012capital}
P.~Beraldi, M.~E. Bruni, and A.~Violi.
\newblock Capital rationing problems under uncertainty and risk.
\newblock \emph{Computational Optimization and Applications}, 51\penalty0 (3):\penalty0 1375--1396, 2012.

\bibitem[Booth et~al.(2016)Booth, Tran, and Beck]{booth2016logic}
K.~E. Booth, T.~T. Tran, and J.~C. Beck.
\newblock Logic-based decomposition methods for the travelling purchaser problem.
\newblock In \emph{Integration of AI and OR Techniques in Constraint Programming: 13th International Conference, CPAIOR 2016, Banff, AB, Canada, May 29-June 1, 2016, Proceedings 13}, pages 55--64. Springer, 2016.

\bibitem[Botton et~al.(2013)Botton, Fortz, Gouveia, and Poss]{botton2013benders}
Q.~Botton, B.~Fortz, L.~Gouveia, and M.~Poss.
\newblock Benders decomposition for the hop-constrained survivable network design problem.
\newblock \emph{INFORMS journal on computing}, 25\penalty0 (1):\penalty0 13--26, 2013.

\bibitem[Broersma et~al.(2013)Broersma, Fomin, van’t Hof, and Paulusma]{broersma2013exact}
H.~Broersma, F.~V. Fomin, P.~van’t Hof, and D.~Paulusma.
\newblock Exact algorithms for finding longest cycles in claw-free graphs.
\newblock \emph{Algorithmica}, 65\penalty0 (1):\penalty0 129--145, 2013.

\bibitem[Church and ReVelle(1974)]{church1974maximal}
R.~Church and C.~ReVelle.
\newblock The maximal covering location problem.
\newblock In \emph{Papers of the regional science association}, volume~32, pages 101--118. Springer-Verlag Berlin/Heidelberg, 1974.

\bibitem[Church et~al.(2004)Church, Scaparra, and Middleton]{Church2004}
R.~L. Church, M.~P. Scaparra, and R.~S. Middleton.
\newblock Identifying critical infrastructure: The median and covering facility interdiction problems.
\newblock \emph{Annals of the Association of American Geographers}, 94\penalty0 (3):\penalty0 491--502, 2004.

\bibitem[Clark et~al.(2006)Clark, Neto, and Toso]{clark2006multi}
A.~R. Clark, R.~M. Neto, and E.~A. Toso.
\newblock Multi-period production setup-sequencing and lot-sizing through atsp subtour elimination and patching.
\newblock In \emph{Proceedings of the 25th workshop of the UK planning and scheduling special interest group. University of Nottingham}, pages 80--87, 2006.

\bibitem[Codato and Fischetti(2006)]{codato2006combinatorial}
G.~Codato and M.~Fischetti.
\newblock Combinatorial benders' cuts for mixed-integer linear programming.
\newblock \emph{Operations Research}, 54\penalty0 (4):\penalty0 756--766, 2006.

\bibitem[Cornu{\'e}jols(2001)]{cornuejols2001combinatorial}
G.~Cornu{\'e}jols.
\newblock \emph{Combinatorial optimization: Packing and covering}.
\newblock SIAM, 2001.

\bibitem[Desrochers and Laporte(1991)]{desrochers1991improvements}
M.~Desrochers and G.~Laporte.
\newblock Improvements and extensions to the miller-tucker-zemlin subtour elimination constraints.
\newblock \emph{Operations Research Letters}, 10\penalty0 (1):\penalty0 27--36, 1991.

\bibitem[Edmonds and Fulkerson(1970)]{edmonds1970bottleneck}
J.~Edmonds and D.~R. Fulkerson.
\newblock Bottleneck extrema.
\newblock \emph{Journal of Combinatorial Theory}, 8\penalty0 (3):\penalty0 299--306, 1970.

\bibitem[Elfar et~al.(2018)Elfar, Talebpour, and Mahmassani]{elfar2018machine}
A.~Elfar, A.~Talebpour, and H.~S. Mahmassani.
\newblock Machine learning approach to short-term traffic congestion prediction in a connected environment.
\newblock \emph{Transportation Research Record}, 2672\penalty0 (45):\penalty0 185--195, 2018.

\bibitem[Ford~Jr and Fulkerson(1956)]{ford1956maximal}
L.~R. Ford~Jr and D.~R. Fulkerson.
\newblock Maximal flow through a network.
\newblock \emph{Canadian journal of Mathematics}, 8:\penalty0 399--404, 1956.

\bibitem[Fouilhoux and Mahjoub(2006)]{fouilhoux2006polyhedral}
P.~Fouilhoux and A.~R. Mahjoub.
\newblock Polyhedral results for the bipartite induced subgraph problem.
\newblock \emph{Discrete Applied Mathematics}, 154\penalty0 (15):\penalty0 2128--2149, 2006.

\bibitem[Gallo et~al.(2009)Gallo, Hammer, and Simeone]{gallo2009quadratic}
G.~Gallo, P.~L. Hammer, and B.~Simeone.
\newblock Quadratic knapsack problems.
\newblock In \emph{Combinatorial optimization}, pages 132--149. Springer, 2009.

\bibitem[Gao et~al.(2023)Gao, Wei, and Walteros]{gao2023exact}
C.~Gao, N.~Wei, and J.~L. Walteros.
\newblock An exact approach for solving pickup-and-delivery traveling salesman problems with neighborhoods.
\newblock \emph{Transportation Science}, 57\penalty0 (6):\penalty0 1560--1580, 2023.

\bibitem[Gao and Kleywegt(2023)]{gao2023distributionally}
R.~Gao and A.~Kleywegt.
\newblock Distributionally robust stochastic optimization with wasserstein distance.
\newblock \emph{Mathematics of Operations Research}, 48\penalty0 (2):\penalty0 603--655, 2023.

\bibitem[G{\"u}nl{\"u}k(1999)]{gunluk1999branch}
O.~G{\"u}nl{\"u}k.
\newblock A branch-and-cut algorithm for capacitated network design problems.
\newblock \emph{Mathematical Programming}, 86:\penalty0 17--39, 1999.

\bibitem[Haynes et~al.(2013)Haynes, Hedetniemi, and Slater]{haynes2013fundamentals}
T.~W. Haynes, S.~Hedetniemi, and P.~Slater.
\newblock \emph{Fundamentals of domination in graphs}.
\newblock CRC press, 2013.

\bibitem[Hooker(2023)]{hooker2023logic}
J.~Hooker.
\newblock \emph{Logic-based benders decomposition: theory and applications}.
\newblock Springer Nature, 2023.

\bibitem[Hooker and Ottosson(2003)]{hooker2003logic}
J.~N. Hooker and G.~Ottosson.
\newblock Logic-based benders decomposition.
\newblock \emph{Mathematical Programming}, 96\penalty0 (1):\penalty0 33--60, 2003.

\bibitem[Israeli and Wood(2002)]{israeli2002shortest}
E.~Israeli and R.~K. Wood.
\newblock Shortest-path network interdiction.
\newblock \emph{Networks: An International Journal}, 40\penalty0 (2):\penalty0 97--111, 2002.

\bibitem[Kara et~al.(2004)Kara, Laporte, and Bektas]{kara2004note}
I.~Kara, G.~Laporte, and T.~Bektas.
\newblock A note on the lifted miller--tucker--zemlin subtour elimination constraints for the capacitated vehicle routing problem.
\newblock \emph{European Journal of Operational Research}, 158\penalty0 (3):\penalty0 793--795, 2004.

\bibitem[Karger(2000)]{karger2000minimum}
D.~R. Karger.
\newblock Minimum cuts in near-linear time.
\newblock \emph{Journal of the ACM (JACM)}, 47\penalty0 (1):\penalty0 46--76, 2000.

\bibitem[Karp(2009)]{karp2009reducibility}
R.~M. Karp.
\newblock Reducibility among combinatorial problems.
\newblock In \emph{50 Years of Integer Programming 1958-2008: from the Early Years to the State-of-the-Art}, pages 219--241. Springer, 2009.

\bibitem[Ke et~al.(2017)Ke, Meng, Finley, Wang, Chen, Ma, Ye, and Liu]{ke2017lightgbm}
G.~Ke, Q.~Meng, T.~Finley, T.~Wang, W.~Chen, W.~Ma, Q.~Ye, and T.-Y. Liu.
\newblock Lightgbm: A highly efficient gradient boosting decision tree.
\newblock \emph{Advances in neural information processing systems}, 30, 2017.

\bibitem[Kheirkhah et~al.(2016)Kheirkhah, Navidi, and Bidgoli]{kheirkhah2016improved}
A.~Kheirkhah, H.~Navidi, and M.~M. Bidgoli.
\newblock An improved benders decomposition algorithm for an arc interdiction vehicle routing problem.
\newblock \emph{IEEE Transactions on Engineering Management}, 63\penalty0 (2):\penalty0 259--273, 2016.

\bibitem[Koopmans and Beckmann(1957)]{koopmans1957assignment}
T.~C. Koopmans and M.~Beckmann.
\newblock Assignment problems and the location of economic activities.
\newblock \emph{Econometrica: journal of the Econometric Society}, pages 53--76, 1957.

\bibitem[Laporte et~al.(1986)Laporte, Nobert, and Arpin]{laporte1986exact}
G.~Laporte, Y.~Nobert, and D.~Arpin.
\newblock An exact algorithm for solving a capacitated location-routing problem.
\newblock \emph{Annals of Operations Research}, 6:\penalty0 291--310, 1986.

\bibitem[Li(2000)]{li2000longest}
M.~Li.
\newblock Longest cycles in almost claw-free graphs.
\newblock \emph{Graphs and Combinatorics}, 16\penalty0 (2):\penalty0 177--198, 2000.

\bibitem[Li et~al.(2022)Li, C{\^o}t{\'e}, Callegari-Coelho, and Wu]{li2022novel}
Y.~Li, J.-F. C{\^o}t{\'e}, L.~Callegari-Coelho, and P.~Wu.
\newblock Novel formulations and logic-based benders decomposition for the integrated parallel machine scheduling and location problem.
\newblock \emph{INFORMS journal on computing}, 34\penalty0 (2):\penalty0 1048--1069, 2022.

\bibitem[Loh(2011)]{loh2011classification}
W.-Y. Loh.
\newblock Classification and regression trees.
\newblock \emph{Wiley interdisciplinary reviews: data mining and knowledge discovery}, 1\penalty0 (1):\penalty0 14--23, 2011.

\bibitem[Lozano and Smith(2017)]{lozano2017backward}
L.~Lozano and J.~C. Smith.
\newblock A backward sampling framework for interdiction problems with fortification.
\newblock \emph{INFORMS Journal on Computing}, 29\penalty0 (1):\penalty0 123--139, 2017.

\bibitem[Luo and Mehrotra(2020)]{luo2020distributionally}
F.~Luo and S.~Mehrotra.
\newblock Distributionally robust optimization with decision dependent ambiguity sets.
\newblock \emph{Optimization Letters}, 14\penalty0 (8):\penalty0 2565--2594, 2020.

\bibitem[Ma et~al.(2015)Ma, Yu, Wang, and Wang]{ma2015large}
X.~Ma, H.~Yu, Y.~Wang, and Y.~Wang.
\newblock Large-scale transportation network congestion evolution prediction using deep learning theory.
\newblock \emph{PloS one}, 10\penalty0 (3):\penalty0 e0119044, 2015.

\bibitem[Maragno et~al.(2025)Maragno, Wiberg, Bertsimas, Birbil, den Hertog, and Fajemisin]{maragno2025mixed}
D.~Maragno, H.~Wiberg, D.~Bertsimas, {\c{S}}.~{\.I}. Birbil, D.~den Hertog, and A.~O. Fajemisin.
\newblock Mixed-integer optimization with constraint learning.
\newblock \emph{Operations Research}, 73\penalty0 (2):\penalty0 1011--1028, 2025.

\bibitem[Marzo et~al.(2022)Marzo, Melo, Ribeiro, and Santos]{marzo2022new}
R.~G. Marzo, R.~A. Melo, C.~C. Ribeiro, and M.~C. Santos.
\newblock New formulations and branch-and-cut procedures for the longest induced path problem.
\newblock \emph{Computers \& Operations Research}, 139:\penalty0 105627, 2022.

\bibitem[McCormick(1976)]{mccormick1976computability}
G.~P. McCormick.
\newblock Computability of global solutions to factorable nonconvex programs: Part i—convex underestimating problems.
\newblock \emph{Mathematical programming}, 10\penalty0 (1):\penalty0 147--175, 1976.

\bibitem[Mohajerin~Esfahani and Kuhn(2018)]{mohajerin2018data}
P.~Mohajerin~Esfahani and D.~Kuhn.
\newblock Data-driven distributionally robust optimization using the wasserstein metric: Performance guarantees and tractable reformulations.
\newblock \emph{Mathematical Programming}, 171\penalty0 (1):\penalty0 115--166, 2018.

\bibitem[Mohamed et~al.(2023)Mohamed, Klibi, Sadykov, {\c{S}}en, and Vanderbeck]{mohamed2023two}
I.~B. Mohamed, W.~Klibi, R.~Sadykov, H.~{\c{S}}en, and F.~Vanderbeck.
\newblock The two-echelon stochastic multi-period capacitated location-routing problem.
\newblock \emph{European journal of operational research}, 306\penalty0 (2):\penalty0 645--667, 2023.

\bibitem[Nobili and Sassano(1989)]{nobili1989facets}
P.~Nobili and A.~Sassano.
\newblock Facets and lifting procedures for the set covering polytope.
\newblock \emph{Mathematical Programming}, 45\penalty0 (1):\penalty0 111--137, 1989.

\bibitem[O'{K}elly(1987)]{o1987quadratic}
M.~E. O'{K}elly.
\newblock A quadratic integer program for the location of interacting hub facilities.
\newblock \emph{European journal of operational research}, 32\penalty0 (3):\penalty0 393--404, 1987.

\bibitem[Oxley(2006)]{oxley2006matroid}
J.~G. Oxley.
\newblock \emph{Matroid theory}, volume~3.
\newblock Oxford University Press, USA, 2006.

\bibitem[Park et~al.(2020)Park, Nielsen, and Moon]{park2020unmanned}
Y.~Park, P.~Nielsen, and I.~Moon.
\newblock Unmanned aerial vehicle set covering problem considering fixed-radius coverage constraint.
\newblock \emph{Computers \& Operations Research}, 119:\penalty0 104936, 2020.

\bibitem[Rahmaniani et~al.(2017)Rahmaniani, Crainic, Gendreau, and Rei]{rahmaniani2017benders}
R.~Rahmaniani, T.~G. Crainic, M.~Gendreau, and W.~Rei.
\newblock The benders decomposition algorithm: A literature review.
\newblock \emph{European Journal of Operational Research}, 259\penalty0 (3):\penalty0 801--817, 2017.

\bibitem[Salmeron et~al.(2004)Salmeron, Wood, and Baldick]{salmeron2004analysis}
J.~Salmeron, K.~Wood, and R.~Baldick.
\newblock Analysis of electric grid security under terrorist threat.
\newblock \emph{IEEE Transactions on power systems}, 19\penalty0 (2):\penalty0 905--912, 2004.

\bibitem[S{\'a}nchez-Garc{\'\i}a et~al.(1998)S{\'a}nchez-Garc{\'\i}a, Sobr{\'o}n, and Vitoriano]{sanchez1998set}
M.~S{\'a}nchez-Garc{\'\i}a, M.~I. Sobr{\'o}n, and B.~Vitoriano.
\newblock On the set covering polytope: Facets with coefficients in $\{$0, 1, 2, 3$\}$.
\newblock \emph{Annals of Operations Research}, 81\penalty0 (0):\penalty0 343--356, 1998.

\bibitem[Sassano(1989)]{sassano1989facial}
A.~Sassano.
\newblock On the facial structure of the set covering polytope.
\newblock \emph{Mathematical Programming}, 44\penalty0 (1-3):\penalty0 181--202, 1989.

\bibitem[Saxena et~al.(2010)Saxena, Goyal, and Lejeune]{saxena2010mip}
A.~Saxena, V.~Goyal, and M.~A. Lejeune.
\newblock Mip reformulations of the probabilistic set covering problem.
\newblock \emph{Mathematical programming}, 121\penalty0 (1):\penalty0 1--31, 2010.

\bibitem[Smith and Song(2020)]{smith2020survey}
J.~C. Smith and Y.~Song.
\newblock A survey of network interdiction models and algorithms.
\newblock \emph{European Journal of Operational Research}, 283\penalty0 (3):\penalty0 797--811, 2020.

\bibitem[Smith et~al.(1977)Smith, Srinivasan, and Thompson]{smith1977computational}
T.~H. Smith, V.~Srinivasan, and G.~Thompson.
\newblock Computational performance of three subtour elimination algorithms for solving asymmetric traveling salesman problems.
\newblock In \emph{Annals of Discrete Mathematics}, volume~1, pages 495--506. Elsevier, 1977.

\bibitem[Song et~al.(2014)Song, Luedtke, and K{\"u}{\c{c}}{\"u}kyavuz]{song2014chance}
Y.~Song, J.~R. Luedtke, and S.~K{\"u}{\c{c}}{\"u}kyavuz.
\newblock Chance-constrained binary packing problems.
\newblock \emph{INFORMS Journal on Computing}, 26\penalty0 (4):\penalty0 735--747, 2014.

\bibitem[Srivastava and Bullo(2014)]{srivastava2014knapsack}
V.~Srivastava and F.~Bullo.
\newblock Knapsack problems with sigmoid utilities: Approximation algorithms via hybrid optimization.
\newblock \emph{European Journal of Operational Research}, 236\penalty0 (2):\penalty0 488--498, 2014.

\bibitem[Tanner and Ntaimo(2010)]{tanner2010iis}
M.~W. Tanner and L.~Ntaimo.
\newblock Iis branch-and-cut for joint chance-constrained stochastic programs and application to optimal vaccine allocation.
\newblock \emph{European Journal of Operational Research}, 207\penalty0 (1):\penalty0 290--296, 2010.

\bibitem[Udell and Boyd(2013)]{udell2013maximizing}
M.~Udell and S.~Boyd.
\newblock Maximizing a sum of sigmoids.
\newblock \emph{Optimization and Engineering}, pages 1--25, 2013.

\bibitem[Wei and Walteros(2022)]{wei2022integer}
N.~Wei and J.~L. Walteros.
\newblock Integer programming methods for solving binary interdiction games.
\newblock \emph{European Journal of Operational Research}, 302\penalty0 (2):\penalty0 456--469, 2022.

\bibitem[Wei and Walteros(2025)]{wei2025supervalid}
N.~Wei and J.~L. Walteros.
\newblock On supervalid inequalities for binary interdiction games.
\newblock \emph{Mathematical Programming}, 212\penalty0 (1):\penalty0 437--478, 2025.

\bibitem[Wei et~al.(2021)Wei, Walteros, and Pajouh]{wei2021integer}
N.~Wei, J.~L. Walteros, and F.~M. Pajouh.
\newblock Integer programming formulations for minimum spanning tree interdiction.
\newblock \emph{INFORMS Journal on Computing}, 33\penalty0 (4):\penalty0 1461--1480, 2021.

\bibitem[Yao et~al.(2007)Yao, Edmunds, Papageorgiou, and Alvarez]{yao2007trilevel}
Y.~Yao, T.~Edmunds, D.~Papageorgiou, and R.~Alvarez.
\newblock Trilevel optimization in power network defense.
\newblock \emph{IEEE Transactions on Systems, Man, and Cybernetics, Part C (Applications and Reviews)}, 37\penalty0 (4):\penalty0 712--718, 2007.

\end{thebibliography}

\newpage
\appendix
\section{Mathematical Proofs}
\label{sec:proofs}



\alg*
\pfstart
\textbf{(A)} Items \ref{alg:a1}, \ref{alg:a3}, \ref{alg:a5} are immediate. For \ref{alg:a4}, if $\Omega$ is upper and $T\in\nc{\Omega}$, then any $T'\subsetneq T$ cannot lie in $\Omega$ (else upperness would force $T\in\Omega$), hence $T'\in\nc{\Omega}$, showing $\nc{\Omega}$ is lower. The reverse direction is symmetric.

\smallskip
\noindent
\textbf{(B)} Items \ref{alg:b1}, \ref{alg:b2}, \ref{alg:b5} are immediate from the definition of $\ec$ as an involutive bijection on elements; \ref{alg:b4} follows since complement reverses inclusion on elements.

\smallskip
\noindent
\textbf{(C)} Items \ref{alg:c1}–\ref{alg:c5} are standard from the definition of cut operator. For \ref{alg:c6}, if $\emptyset\in\Omega$ then $\mathcal C(\Omega)=\emptyset$ and $\mathcal C\mathcal C(\Omega)=\pset(\Delta)=\uparrow\Omega$. Otherwise, $U\in\mathcal C\mathcal C(\Omega)$ iff $U$ intersects every $S$ that intersects every $T\in\Omega$, which is equivalent to $U\supseteq T$ for some $T\in\Omega$, i.e., $U\in\uparrow\Omega$.

\smallskip
\noindent
\textbf{(E)} Items \ref{alg:e1}--\ref{alg:e4} follow from Lemma~\ref{lem:contra}. For \ref{alg:e5}, $T\in\nc{\uparrow\Omega}\Rightarrow T\notin\uparrow\Omega \Rightarrow T\notin\Omega \Rightarrow T\in\nc{\Omega}\Rightarrow T\in\downarrow\nc{\Omega}$. The second inclusion is analogous. For \ref{alg:e6}, $T\in m\nc{(\Omega)}$ implies $T\notin M(\Omega)$, hence $T\in\nc{M(\Omega)}$ (and symmetrically for the other inclusion). For \ref{alg:e7}, if $T\in\mathcal G\nc{(\Omega)}$, then for all $T'\notin\Omega$ we have $T\cup T'\neq\Delta$, which implies $\Delta\setminus T\in\Omega$, hence $T\notin \mathcal C(\Omega)$ and $\mathcal G\nc{(\Omega)}\subseteq\nc{\mathcal C(\Omega)}$. For the equality part, we first show sufficiency. Take any $T \notin \mathcal C(\Omega)$, there exists some $T' \in \Omega$ such that $T' \cap T = \emptyset$. Suppose $T \notin \mathcal G(\nc\Omega)$, then there exists some $T'' \notin \Omega$ such that $T \cup T''=\Delta$. In particular, we have $T' \subseteq T''$. But, $T' \in \Omega$ and $T'' \notin \Omega$ contradicts that $\Omega$ is an upper system. For the necessity, we prove the contrapositive: suppose $\Omega$ is not an upper system, then $\nc{\mathcal C(\Omega)} \cap \nc{\mathcal G\nc(\Omega)} \neq \emptyset$. Note that $\Omega$ is not an upper system implies that there exist $T \in \Omega$ and $T' \notin \Omega$ such that $T \subseteq T'$. We show that every $T''$ sandwiched between $T$ and $T'$ entails that its complement $\Delta \setminus T''$ belongs to the intersection $\nc{\mathcal C(\Omega)} \cap \nc{\mathcal G\nc(\Omega)}$ by the following
  $$T \subseteq T'' \Longleftrightarrow T \subseteq \Delta\setminus (\Delta \setminus T'') \Longrightarrow T \cap (\Delta \setminus T'') = \emptyset \Longrightarrow \Delta \setminus T'' \notin \mathcal C(\Omega),$$
  $$T' \supseteq T'' \Longleftrightarrow T' \supseteq \Delta\setminus (\Delta \setminus T'') \Longrightarrow T' \cup (\Delta \setminus T'') = \Delta \Longrightarrow \Delta \setminus T'' \notin \mathcal G\nc(\Omega),$$
which completes the proof of necessity. The proof for $\mathcal C\nc (\Omega) \subseteq \nc\mathcal G(\Omega)$ is symmetric.

\smallskip
\noindent
\textbf{(D)} Follows from the identities in (E) and the properties in (C) by duality. For instance, \ref{alg:d1}–\ref{alg:d2} follow from \ref{alg:c1}–\ref{alg:c2} via $\mathcal G=\ec{\mathcal C\ec}$; \ref{alg:d3} from \ref{alg:c3}; \ref{alg:d4} from \ref{alg:c4} together with $\ec{m\ec}=M$; \ref{alg:d5} from \ref{alg:c5} using that $\ec$ is order-preserving on set systems; and \ref{alg:d6} from \ref{alg:c6} plus $\uparrow\ec=\ec{\downarrow}$.
\pfend

\muembd*
\pfstart
For the first inclusion, take any $T \in \mathcal C\ec{\nc(\Omega)}$ and suppose $T \notin \Omega$. By definition, we have $T \in \nc \Omega$ and $\Delta \setminus T \in \ec{\nc \Omega}$. On the other hand, $T \in \mathcal C\ec{\nc(\Omega)}$ means $T \cap S \neq \emptyset$ for all $S \in \ec{\nc\Omega}$. In particular, $T \cap (\Delta \setminus T) \neq \emptyset$, a contradiction. 
To prove $\mathcal C\ec{\nc(\Omega)}$ is the largest inner approximation, we show that any upper system $\Omega' \subseteq \Omega$ is a subset of $\mathcal C\ec{\nc(\Omega)}$. Suppose otherwise, take $T \in \Omega' \subseteq \Omega$ but $T \notin \mathcal C\ec{\nc(\Omega)}$. The latter implies that for some $S \in \ec{\nc\Omega}$, we have $T \subseteq \Delta \setminus S \in \nc\Omega$, where the membership also implies $\Delta \setminus S \notin \Omega$. However, since $\Omega'$ is an upper system, we have $\Delta \setminus S \in \Omega' \subseteq \Omega$, a contradiction.
For the second inclusion, $\uparrow \Omega$ is clearly the tightest outer approximation of $\Omega$ by definition. Applying the above argument with the system $\uparrow \Omega$, we obtain
$$\mathcal C\ec{\nc{\uparrow(\Omega)}} = \uparrow\Omega,$$
where the equality is due to $\uparrow\Omega$ is upper-closed.
Finally, since the systems on both sides are upper-closed, the equality holds on both sides if and only if $\Omega$ is also upper-closed.
\pfend

\embdvar*
\pfstart
Replacing $\Omega$ in the first inclusion of Theorem~\ref{thm:muembd} with each of the right-side sets, we can derive these claims using the commutativity rule \ref{alg:e1} and the self-inverse rules \ref{alg:a5} and \ref{alg:b5}. Moreover, they are the tightest upper embeddings by Theorem \ref{thm:muembd}.
\pfend

\ce*
\pfstart
By definition, $S \in \mathcal C(\Omega)$ if and only if $S \cap T \neq \emptyset$ for all $T \in \Omega$. Thus, the vector representation $x_S$ is feasible if and only if it satisfies all the covering inequalities. On the other hand, we have $S \in \ec{\mathcal C(\Omega)}$ if and only if $\Delta \setminus S \in \mathcal C(\Omega)$. Hence, vector $y_S:=1-x_S$ indicates the complement structure $\Delta \setminus S$ and intersects every $T \in \Omega$. Substituting $x=1-y$ in the covering inequalities produces the elimination constraints. In particular, this also proves the last claim.
\pfend

\reform*
\pfstart
These approximation results directly follow Theorem~\ref{thm:muembd}, Corollary~\ref{coro:lembd}, and Lemma~\ref{lem:ce}. The value relationships follow the definition of inner and outer approximations.
\pfend

\sep*
\pfstart
Suppose \(T\in\Omega\). Then no violated inequality is needed. Otherwise, the separation algorithm constructs an extremal infeasible set from \(T\), namely a maximal (or minimal ) set in \(\nc\Omega\) when \(\Omega\) is upper-closed (or lower-closed). This can be done by a linear search over at most \(|\Delta|\) elements, using one oracle call per element. The resulting extremal set indexes a violated covering or elimination inequality, which proves the claimed complexity.
\pfend

\faceti*
\pfstart
For every fixed index subset $I \subseteq [n]$, the transformation $\theta_I(\cdot)$ is clearly affine and injective by definition. To show $\theta_I$ is also surjective, we represent an arbitrary $x' \in \conv(\mathcal X_{\theta_I})$ as the convex combination $\sum_{j \in J}\lambda_j x'_j$ for some extreme points $x'_j$ of $\conv(\mathcal X_{\theta_I})$ and define $x = \sum_{j \in J}\lambda_j \theta_I(x'_j)$. Clearly, we have $\theta_I(x) = x'$ since $\theta_I$ is affine and is an involution. It remains to show $x \in \mathcal X$. By construction, each $x'_j$ is binary and satisfies $g(\theta_I(x'_j)) \leq 0$, implying $\theta_I(x'_j) \in \mathcal X$. Hence, the convex combination $x$ belongs to the convex hull $\conv(\mathcal X)$. Since affine bijections preserve faces and their dimensions, the statement follows.
\pfend


\facetthm*
\pfstart
This statement is a paraphrase of the Statements $1$ and $5$ in \cite{balas1989set} (from the matrix perspective) and Proposition 3.4 in \cite{wei2022integer} (from the set system perspective). According to the latter, given an upper system $\Omega$ such that $|T| \geq 2$ for every $T \in m\ec{\nc{(\Omega)}}$, \eqref{eq:pformub01} is facet-defining in the associated set cover problem if and only if for every $a \in \Delta \setminus T$ there exists some $a' \in T$ such that $T \setminus \{a'\} \cup \{a\} \notin \ec{\nc\Omega}$, that is, $\Delta \setminus (T \setminus \{a'\} \cup \{a\}) = (\Delta\setminus T) \setminus \{a\} \cup \{a'\} \in \Omega$. By definition, this means $\Delta \setminus T$ is quasi-feasible. 
\pfend

\func*
\pfstart
For every linear function $\iprod{c,x}$ with index set $[n]$, defining $I:=\{i \in [n] \mid c_i \geq 0\}$ induces the required index bipartition. The described bilinear functions with index set $I \cup J$ can be written as $\iprod{x_{I_1}, R_{I_1 J_1}y_{J_1}} + \iprod{x_{I_2}, R_{I_2 J_2}y_{J_2}}$. Since $R_{I_1 J_1} \geq 0$ and $R_{I_2 J_2}\leq 0$ by assumption, the index bipartition $(I_1 \cup J_1, I_2 \cup J_2)$ satisfies the bimonotone requirement. For submodular functions, we have
$$f(\{i\})-f(\emptyset) \geq f(T \cup \{i\}) - f(T) \geq f(\Delta) - f(\Delta \setminus \{i\})$$
for every $i \in \Delta$ and every $T \subseteq \Delta \setminus \{i\}$. Hence, defining $I:=\{ i \mid f(\Delta) - f(\Delta \setminus \{i\}) \geq 0\}$ induces the required index bipartition. The supermodular case has all above inequality reversed, thus the index bipartition induced by $I:=\{ i \mid f(\{i\}) - f(\emptyset) \geq 0\}$ proves the claim.
\pfend


\bimono*
\pfstart
Let \(x'=\theta(x)\), where \(x'_I=x_I\) and \(x'_J=1-x_J\). Under this change of variables, the transformed feasible set is \(\Omega_\theta\), and the objective \(f(\theta^{-1}(x'))\) is increasing in \(x'\). Therefore, minimizing over \(\Omega_\theta\) is equivalent in optimal value to minimizing over its upper closure \(\uparrow\Omega_\theta\): every point in \(\uparrow\Omega_\theta\setminus\Omega_\theta\) contains some point in \(\Omega_\theta\) that is componentwise no larger, and hence has no larger objective value.
By Theorem~\ref{thm:reform}, the upper-closed system \(\uparrow\Omega_\theta\) admits the covering formulation indexed by
\(
m\ec{\nc(\uparrow\Omega_\theta)}.
\)
Applying the inverse transformation \(x'=\theta(x)\) gives the constraints in \eqref{eq:bimono}. Hence the reformulation is exact. The facet condition follows by applying Proposition~\ref{prop:facet} to the upper system \(\uparrow\Omega_\theta\), together with the facet-preserving property of the flipping map in Lemma~\ref{lem:faceti}.
\pfend


\interval*
\pfstart
It is immediate that $\Omega \subseteq \uparrow\Omega \cap \downarrow\Omega$.  
To show that the intersection forms an interval system, take $T,T' \in \uparrow\Omega \cap \downarrow\Omega$ and any $T''$ with $T \subseteq T'' \subseteq T'$.  
Since $T \in \uparrow\Omega$, there exists $S \in \Omega$ with $S \subseteq T \subseteq T''$, hence $T'' \in \uparrow\Omega$.  
Similarly, because $T' \in \downarrow\Omega$, there exists $S' \in \Omega$ with $T' \subseteq S'$, and thus $T'' \subseteq S'$, so $T'' \in \downarrow\Omega$.  
Therefore $T'' \in \uparrow\Omega \cap \downarrow\Omega$, proving closure under intervals.  
For minimality, let $\Omega'$ be any interval system with $\Omega \subseteq \Omega'$.  
Take any $T \in \uparrow\Omega \cap \downarrow\Omega$. Then there exist $T_1,T_2 \in \Omega$ with $T_1 \subseteq T \subseteq T_2$.  
Since $\Omega \subseteq \Omega'$ and $\Omega'$ is interval, it follows that $T \in \Omega'$.  
Thus $\Omega' \supseteq \uparrow\Omega \cap \downarrow\Omega$, showing minimality.  
For the second statement, if $\Omega$ is interval, then by minimality we must have $\Omega \supseteq \uparrow\Omega \cap \downarrow\Omega$, and the reverse inclusion holds trivially.
Conversely, if $\Omega = \uparrow\Omega \cap \downarrow\Omega$, then $\Omega$ equals an interval system, hence itself is interval.  
\pfend

\end{document}